\documentclass[12pt]{amsart} 
\usepackage[mathscr]{eucal} 
\usepackage{amsmath,amsfonts} 
\usepackage{amssymb}
\parskip=\smallskipamount 
\newtheorem{theorem}{Theorem}[section] 
\newtheorem{proposition}[theorem]{Proposition} 
\newtheorem{corollary}[theorem]{Corollary} 
\newtheorem{lemma}[theorem]{Lemma} 
\theoremstyle{definition} 
 
\newtheorem{example}[theorem]{Example}

\newtheorem{remark}[theorem]{Remark} 
\newtheorem{remarks}[theorem]{Remarks}
\makeatletter 
\@addtoreset{equation}{section} 
\makeatother

\makeatletter
\@namedef{subjclassname@2020}{\textup{2020} Mathematics Subject Classification}
\makeatother

\newcommand{\CC}{{\mathbb C}} 
\newcommand{\NN}{{\mathbb N}}

\newcommand{\RR}{{\mathbb R}} 
 
\newcommand{\cA}{{\mathcal A}} 
\newcommand{\cB}{{\mathcal B}} 
\newcommand{\cC}{{\mathcal C}}

\newcommand{\Ra}{\Rightarrow}

\newcommand{\ra}{\rightarrow} 
\newcommand{\ol}{\overline}

\let\phi=\varphi

\newcommand{\tl}{\tilde}

\begin{document} 
\title[Isometric Representations]{Isometric Representations of Calibrated Ordered Spaces on $C(X)$}\thanks{}
 
 \date{\today}
 
 \author[S. Ay]{Serdar Ay}
 \address{Department of Mathematics, At{\.i}l{\.i}m University, 06830 {\.I}ncek, Ankara, 
Turkey}
 \email{serdar.ay@atilim.edu.tr}
 

\begin{abstract} 
The problem of characterizing normed ordered spaces which admit a representation in the algebraic, order and 
norm sense as a subspace of $C(X)$, the space of all continuous functions on a compact Hausdorff space 
is a classical problem that has been considered by many authors. In this article 
we consider the more general case of calibrated ordered spaces, that is, 
ordered spaces with a specified family of seminorms generating its topology. For such spaces 
equivalent conditions on representability as a subspace of $C(X)$ for some 
locally compact Hausdorff space $X$, in the algebraic, order and seminorm sense are stated and proved. 
Some characterizations appear to be new even in the normed case. 
As an application 
of the main theorems, we state and prove a characterization of norm additivity property 
of two positive functionals. 
\end{abstract} 

\maketitle 
 
\section*{Introduction}

The problem of characterizing normed ordered spaces which admit a representation in the algebraic, order and 
norm sense as a subspace of $C(X)$, the space of all continuous functions on a compact Hausdorff space 
is a classical problem that has been considered by many authors. In the article \cite{Kakutani}, isometrically isomorphic 
representations of AM (Abstract M)-spaces are obtained. In \cite{Nachbin} and \cite{Riedl} characterizations 
in terms of semi positive and semi negative elements and fullness of the unit ball for 
isometrically isomorphic representability on $C(X)$, where $X$ is a compact Hausdorff space, were made respectively. 
In \cite{Kadison} the case of ordered normed spaces with an order unit was 
considered and isometrically isomorphically representability of the self adjoint part of a unital $C^*$-algebra 
as a space of continuous functions $C(X)$ was shown. For complex normed order unit spaces see \cite{PaulsenTomforde}, 
see also \cite{Dosi} and \cite{Asadi} for the related notion of a local order unit space. 
Isometrically isomorphic representability as a subspace of a normed Riesz space, and normed M-Riesz space 
is discussed in \cite{vanGaans}. In the recent article \cite{deBruyn1} characterizations of 
injective, positive and continuous representability of a real ordered topological vector space $L$ on the space 
$C(X)$ where $X$ is locally compact are obtained. 
In \cite{Schaefer1}, 
\cite{Schaefer2}, the closely related result of representing a normal ordered space in algebraic, order 
and topological sense as a subspace of $C(X)$, 
where $X$ is a locally compact Hausdorff space 
was obtained, see also \cite{Aliprantis} and \cite{Riedl}. Probably there are many other sources treating the subject. 

In this article we consider the case of calibrated ordered spaces, that is, 
ordered spaces with a specified family of seminorms generating its topology, 
see sections \ref{s:nps} and \ref{s:isoschrep} for precise definitions. For such spaces 
equivalent conditions on representability as a subspace of $C(X)$ for some 
locally compact Hausdorff space $X$, in the algebraic, order and seminorm sense are stated and proved. 
As an application 
of the main theorems, we state and prove a characterization of seminorm additivity property 
of two positive functionals, see Theorem \ref{t:flsnorm}. 

It seems that isometrically isomorphically representability 
of calibrated ordered spaces on some $C(X)$ was not considered till now. It can be considered as a natural 
continuation of Schaefer's construction on representations of normal ordered spaces, e.g. see the proof 
of Theorems \ref{t:posext}, \ref{t:isorepcx}, \ref{t:isorepsnfl}, \ref{t:isorepstr} and \ref{t:isorepquot}. In this 
setting we obtain generalizations of 
known statements on characterizations of isometrically isomorphically representability of normed ordered spaces 
as well as some characterizations which appear to be new 
even in the normed case, see the main theorems \ref{t:isorepcx}, \ref{t:isorepsnfl}, 
\ref{t:isorepstr} and \ref{t:isorepquot} as well as Remark \ref{r:refofisorep}. An additional 
impetus to write this article was given to us by the recent article \cite{deBruyn1}, where ordered 
topological vector spaces were considered and characterizations of representability in a weaker sense as a subspace of $C(X)$ 
where $X$ is a locally compact Hausdorff space were obtained. This article could be considered 
as a continuation of that article in the locally convex direction. 

In \cite{Asadi} representations of locally ordered $*$-spaces by local order isomorphisms as 
a subspace of $C(X)$, where $X$ is a compactly generated completely Hausdorff topological space, was obtained. The notion 
of local order and in particular a local order unit space was first defined and studied in \cite{Dosi}. 
In this article we give a less restrictive definition of a real local order unit space 
and consequently obtain a more general statement concerning their representations, see 
Section \ref{s:nps}, Theorem \ref{t:isorepstr} and Remark \ref{r:refofisorep}.

On the other hand, the notion of VE (Vector Euclidean) and VH (Vector Hilbert) spaces, introduced by Loynes 
in \cite{Loynes1}, \cite{Loynes2} and \cite{Loynes3} are closely related 
to complex ordered spaces with a locally convex topology. 
These are 
vector spaces equipped with an ordered space valued inner product and 
therefore generalizations of locally Hilbert $C^*$-modules 
and of Hilbert $C^*$-modules, in particular of those which 
are over $C(X)$ for $X$ a locally compact Hausdorff space. See also the articles  
\cite{Ciurdariu}, \cite{GasparGaspar}, 
\cite{GasparGaspar2}, \cite{GasparGaspar1}, \cite{Saworotnow}, \cite{Weron}, \cite{WeronChobanyan} as well as the more recent \cite{Ay1}, 
\cite{AyGheondea1}, 
\cite{AyGheondea2}, \cite{AyGheondea3}, \cite{Gheondea}, \cite{Gheondeaugur} on this subject. We think that an investigation of 
representations of calibrated ordered spaces on some $C(X)$ might be helpful in 
understanding the structures of VE and VH-spaces better and in particular 
bridging the gap between them and locally Hilbert $C^*$-modules.   

Central to our proofs are Schaefer's construction, 
which we state with a clause on isometric representations of positive elements, see Section \ref{s:isoschrep}, 
and the Bauer-Nachbin-Namioka Theorem and its corollaries, see Section \ref{s:phb}. To the best of our knowledge, 
this approach was not used before, and in addition to new results, it gives rise to new proofs of some of the known statements.  

We briefly describe the contents of this article. In Section \ref{s:nps} we review 
the definition and some properties of ordered spaces and their functionals. Section \ref{s:phb} reviews 
the statement and a proof of 
a version of the Bauer-Namioka Theorem, which we refer to as the Bauer-Nachbin-Namioka Theorem 
for historical reasons explained in the paragraph preceeding the theorem, see Theorem \ref{t:baunami}. 
In the rest of the section corollaries that will be relevant are stated and proved. 

We define calibrated ordered spaces in Section \ref{s:isoschrep}. We then present and prove a 
version of Schaefer's representation theorem for calibrated normal ordered spaces, with 
a slightly different construction from the one given in \cite{Schaefer1}, see 
Theorem \ref{t:schaf}. The construction in that theorem is used in the main theorems. 

In Section \ref{s:cirpolcs} we state and prove the main theorems of the paper, see 
theorems \ref{t:isorepcx}, \ref{t:isorepsnfl}, \ref{t:isorepstr} and \ref{t:isorepquot}. We first give definitions of some notions 
in the calibrated ordered setting; e.g. locally Riesz spaces. We then 
do preparation for the proof of the main theorems; Theorem \ref{t:posext},   
the discussion after it and Theorem \ref{t:rmlocal} are in this direction. After 
the main theorems we present some examples. Finally, in Section \ref{s:naddpos} we state 
and prove an application of the main theorems characterizing 
pairs of positive functionals with the norm additive property, 
see Theorem \ref{t:flsnorm}. 

\section{Notation and Preliminary Results}\label{s:nps}

\subsection{Ordered Spaces.}\label{ss:as}
A real vector space $E$ is called an \emph{$\mathbb{R}$-ordered space}, or \emph{ordered space}, see e.g. \cite{Aliprantis} if 
$E$ is an $\mathbb{R}$-vector space which has a specified cone $E_+$, i.e. for all $x_1,x_2\in E_+$ and $\alpha,\beta\geq 0$ 
we have $\alpha x_1 + \beta x_2 \in E_+$ so that $E_+$ is a wedge and the wedge is \emph{strict}, i.e. $E_+\cap -E_+ =\{0\}$. For 
complex versions of ordered spaces, called ordered $*$-spaces, see e.g. \cite{PaulsenTomforde} or 
\cite{AyGheondea2}.  

Clearly, any vector subspace $W$ of an ordered space $Z$ can be considered as an ordered subspace 
by considering the cone $W_+:=W\cap E_+$. 

The cone $E_+$ is called 
\emph{generating} if for any element $x\in E$ there exists $y_1, y_2 \in E_+$ such that 
$z=y_1 - y_2$. In this article we do not assume generating property of the cone, unless otherwise is specified. 

The specified cone induces a partial ordering $E$, defined by $x\geq y \iff (x - y) \in E_+ $. 

An $\RR$-linear functional $f\colon E\ra \RR$ is called \emph{positive} if $f(x)\geq 0$ for all $x\in E_+$. 
For arbitrary linear functionals $f,g$ on $E$, $f\geq g$ if $f-g$ is a positive functional on $E$. On the other 
hand, a linear map $\phi\colon E \ra F$ where $E$ and 
$F$ are two ordered spaces is called 
\emph{bipositive} if $x \in E_+ \iff  \phi(x)\in F_+$. Notice that, a bipositive map $\phi$ 
is order preserving, that is, $x \geq y \iff \phi(x) \geq \phi(y)$ 
for $x,y\in E$. 

An element $e\in E$ is called an \emph{order unit} if for any $x\in E$ there exists $\lambda > 0$ such that 
$x\leq \lambda e$. See \cite{Aliprantis}, \cite{Riedl} as well as \cite{PaulsenTomforde}.  

We say that $E$ is \emph {Archimedean} if for any $y\in E$ and $x\in E_+$, 
$y+\lambda x \geq 0$ for all $\lambda > 0$ implies that $y\geq 0$, see 
e.g. \cite{Aliprantis}, \cite{PaulsenTomforde} as well as \cite{AyGheondea4}. 
$E$ is called \emph{almost Archimedean} if for any points $x,y\in E$, the relation 
$-\lambda x \leq y \leq \lambda x$ for all $\lambda > 0$ implies $y=0$, see \cite{Riedl}, p.105. 
Since the cone $E_+$ is strict, it follows that an Archimedean ordered space is almost Archimedean. 

The following proposition can be found on p.105 of \cite{Riedl}. 

\begin{proposition}
Let $E$ be an ordered space and assume that $e\in E$ is an order unit. Then $p_{e}\colon E\ra \RR$ 
given by $p_e(x):=\mathrm{inf} \{\lambda >0 \mid -\lambda e \leq x \leq \lambda e  \} $ is a 
seminorm on $E$, called the \emph{order seminorm} associated to $e$. Moreover, 
$p_e$ is a norm if and only if $E$ is almost Archimedean.
\end{proposition}

Let $E$ be an ordered space and $p$ be a seminorm on it. Then by $E^{'}_{p}$ we denote 
the vector space of all functionals $f\colon E\ra \RR$ for which 
there exists a constant $C\geq 0$ such that $\vert f(x) \vert \leq C p(x)$ for all $x\in E$. 
The infimum of all such $C\geq 0$ is denoted by $\vert f \vert_{p}$ and it is easy to check 
that it is a norm on $E^{'}_{p}$. By the symbol $E^{'}_{p,+}$ we denote the subset of 
all positive functionals of $E^{'}_{p}$. 

For an ordered space $E$ with a locally convex topology, the cone $E_+$ is 
said to be \emph{normal} if there exists a family $\{ N_j \}_{j\in J}$ of neighbourhoods 
of $0$, linearly generating the topology, such that $0\leq y \leq x$ with 
$x\in N_j$ implies $y\in N_j$. Equivalently, there exists a family $\{ p_j \}_{j\in J}$ 
of \emph{increasing seminorms} generating the topology of $E$, that is, 
$p_j(x)\geq p_j(y)$ for any $j\in J$ and any elements $0\leq y \leq x$. See \cite{Aliprantis} or \cite{Schaefer1}.  

In the recent article \cite{Dosi}, the notion of locally ordered $*$-vector spaces was introduced and studied 
in connection with the problem of characterizing some completely positive maps between Archimedean order unit spaces 
and related results. We will consider real locally ordered spaces with a less restrictive 
definition and show that a real version of 
Theorem 2.13 in \cite{Asadi} falls under our Theorem \ref{t:isorepstr}. Below we will provide 
some preliminary definitions and facts on locally ordered spaces; for more details see \cite{Dosi} 
and \cite{Asadi}. 

Consider an ordered space $(E,E_+)$ where 
a family of wedges $\{ C_{\alpha} \}_{\alpha\in A}$ is specified with the condition that $\cap_{\alpha\in A} C_{\alpha} = E_+$. 
Then the pair $(E,\{ C_{\alpha} \}_{\alpha\in A})$ is called a \emph{locally ordered space}. Notice 
that one has a corresponding family of relations $\{\leq_{\alpha}\}_{\alpha\in A}$ defined by $x\leq_{\alpha} y$ if and only if 
$y-x\in C_{\alpha}$. Since the sets $C_{\alpha}$ are not assumed to be strict, in general the relations 
$\leq_{\alpha}$ are reflexive and transitive, but not antisymmetric, thus are \emph{preorders}. Note also 
that if $x\leq y$ for any $x,y\in E$, then we have $x\leq_{\alpha} y$ for all $\alpha \in A$. 

If an element $e\in E$ 
has the property that for any $\alpha\in A$ and $x\in E$ there exists a positive number $r_{\alpha}$ such that 
$r_{\alpha}e-x \in C_{\alpha}$, then $e$ is called a \emph{local order unit}. A locally ordered space $(E,\{ C_{\alpha} \}_{\alpha\in A})$ 
is called 
\emph{L-Archimedean} if for any $\alpha\in A$, $y\in E$ and $x\in C_{\alpha}$, 
$y+\lambda x \geq_{\alpha} 0$ for all $\lambda > 0$ implies that $y \in C_{\alpha}$. In the case when $(E,\{ C_{\alpha} \}_{\alpha\in A})$ 
has a local order unit $e$, an easy argument shows that the space is L-Archimedean if and only if 
for any $\alpha\in A$ and $x\in E$ the relation $\lambda e + x \geq_{\alpha} 0$ for all $\lambda > 0$ implies $x\in C_{\alpha}$. 
Moreover it is easy to see that an L-Archimedean locally ordered space is in particular an Archimedean ordered space. 

\begin{remark}\label{r:dosidef}
We remark that in the definitions in 
\cite{Dosi} and \cite{Asadi}, the specified family of wedges $\{ C_{\alpha} \}_{\alpha\in A}$ must have 
\begin{itemize}
\item[(i)] The family $\{ C_{\alpha} \}_{\alpha\in A}$ is \emph{downward filtered}; i.e. for any 
$\alpha, \beta \in A$ there exists $\gamma\in A$ such that $C_{\gamma}\subseteq C_{\alpha},\,C_{\beta}$. 

\item[(ii)] We have $\cap_{\alpha\in A} C_{\alpha} = E_+$. 
\end{itemize}

Therefore we skip the downward filtered condition and have a less restrictive definition of local 
order unit spaces. We will see that, this definition is sufficient for the purposes of this article, 
see Example \ref{r:cxlocou} and Theorem \ref{t:isorepstr}. 
\end{remark}

In the article we will state a result concerning locally $C^*$-algebras, 
therefore we recall their definition and basic properties as well. A $*$-algebra $\cA$ with a complete 
Hausdorff topology given by a family of
\emph{$C^*$-seminorms}, that is, seminorms $p$ on $\cA$ satisfying the
\emph{$C^*$-condition} $p(a^*a) = p(a)^2$ for all $a\in\cA$, is called a 
\emph{locally $C^*$-algebra} \cite{Inoue} 
(equivalent names are \emph{(Locally Multiplicatively Convex) $LMC^*$-algebras}
\cite{Schmudgen}, \cite{Mallios}, or \emph{$b^*$-algebra}
\cite{Allan}, \cite{Apostol}, 
or \emph{pro $C^*$-algebra} \cite{Voiculescu}, \cite{Phillips}). Any 
$C^*$-seminorm is automatically \emph{submultiplicative},  
$p(ab)\leq p(a)p(b)$ for all $a,b\in\cA$, see \cite{Sebestyen},
and is a $*$-seminorm, $p(a^*)=p(a)$ for all $a\in\cA$. 
Denote the collection of all continuous $C^*$-seminorms 
by $S^*(\cA)$. Then $S^*(\cA)$ is a 
directed set under pointwise maximum seminorm, namely, 
given $p,q\in S^*(\cA)$, letting 
$r(a):=\mbox{max}\{p(a),q(a)\}$ for all $a\in\cA$, then 
$r$ is a continuous $C^*$-seminorm and $p,q\leq r$. 
Locally $C^*$-algebras were studied in \cite{Allan}, \cite{Apostol}, \cite{Gheondea2}, 
\cite{Inoue}, \cite{Phillips}, \cite{Schmudgen} and \cite{Zhuraev}, to cite a
few. 

A locally $C^*$-algebra $\cA$ has a natural cone of positive elements, 
$\cA_+:=\{a^*a \mid a\in \cA  \}$. An element $a\in \cA$ is called 
\emph{self adjoint} or \emph{Hermitian} if $a=a^*$. Considering the real 
vector space $\cA^{h}$ of \emph{Hermitian} 
elements of $\cA$, it is well known that $(\cA^{h},\cA_+)$ is an ordered space, 
see e.g. \cite{AyGheondea2}. A $\CC$ linear functional $f\colon \cA \to \CC$ 
is called positive if $f(a^*a)\geq 0$ for all $a\in \cA$. In particular, given 
$p\in S^*(\cA)$, if a positive functional $f$ has $\vert f \vert_p = 1$, where 
$\vert f \vert_p$ has a similar definition to the case of ordered spaces, 
then we call $f$ a \emph{state} of $\cA$.  


\section{Bauer-Nachbin-Namioka Theorem and its Consequences}\label{s:phb}

In this section we will use some standart results on locally convex spaces; i.e. see \cite{Narici} or 
\cite{Conway}. 

In the year 1957, a necessary and sufficient condition 
for the existence of positive linear extensions of positive linear functionals on real ordered vector spaces was discovered by 
H. Bauer in \cite{Bauer1} and 
also in \cite{Bauer2} and independently by I. Namioka in \cite{Namioka}. Verious versions of this result are known 
today as the Bauer-Namioka Theorem; e.g. for a purely algebraic version see Theorem 1.35 in \cite{Aliprantis},
for continuous extensions on topological vector spaces see 5.4. in \cite{Schaefer1}, for a 
version with extensions majorized by a specified seminorm, see \cite{Bauer2}, Hauptsatz 1. See also \cite{Hustad} and 
the recent article \cite{deBruyn2} on extendability of positive continuous functionals.  

On the other hand, in the year 1949, in an article by L. Nachbin, see \cite{Nachbin}, a necessary and sufficient 
condition for the existence of bounded and positive extensions of a linear functional on a subspace 
was provided, see Theorem 2 of that paper. 
We include this condition as equivalent condition $(4)$ in the statement of the theorem 
below, and thus we refer to this version of the theorem as the \emph{Bauer-Nachbin-Namioka Theorem}. 
This version actually follows as a special case of Theorem 10.12b of \cite{Riedl} except 
for condition $(5)$, but 
for convenience of the reader, we present a proof here. 

\begin{theorem}(Bauer-Nachbin-Namioka Theorem)\label{t:baunami}
Let $E$ be a real ordered vector space and $X$ be a vector subspace of $E$. 
Let $p$ be a seminorm on $E$ and $f\colon X \ra \RR$ be a linear functional. Then TFAE: 

\begin{itemize}
\item[(1)] The functional $f$ has a positive linear extension $\tl{f}\colon E\ra \RR$ with 
$\vert \tl{f}(x) \vert \leq p(x)$ for all $x\in E$. 
\item[(2)] The condition 
\begin{equation*}
x \in X\cap(A+E^+) \Ra f(x)\geq -1
\end{equation*}
holds, where $A:=\{ x\in E \mid p(x)\leq 1 \}$.
\item[(3)] The condition 
\begin{equation*}
x \in X\cap(A-E^+) \Ra f(x)\leq 1
\end{equation*}
holds, where $A:=\{ x\in E \mid p(x)\leq 1 \}$.
\item[(4)] For any $x\in X$ and $y\in E$ with $x\leq y$ the condition $f(x) \leq p(y)$ holds.  
\item[(5)] For any $x\in X$ and $y\in E$ with $y\leq x$ the condition $f(-x) \leq p(y)$ holds.
\end{itemize}

If any, hence all of the equivalent conditions hold, 
then $f\in X^{'}_{p,+}$, with notation as in Section \ref{s:nps}. Moreover, for a nonzero functional 
$g\in X^{'}_{p,+}$, an extension $\tl{g}\in E^{'}_{p,+}$ with 
$\vert g \vert_{p} = \vert \tl{g} \vert_{p}$ exists if and only if 
any, hence all conditions are satisfied for the functional $g/\vert g \vert_{p}$.   
\end{theorem}

\begin{proof}
``$(1)\Ra(2)$": Assume that there exists an extension $\tl{f}$ with 
the specified properties. Note that $X\cap (A+E^+)\neq \emptyset$, since $0$ is always there. Let $x\in X\cap (A+E^+)$. 
Then there are $a\in A$ and $l\in E^+$ with $x=a+l$ and 
\begin{equation*}
f(x)=f(a+l)=\tl{f}(a+l)=\tl{f}(a)+\tl{f}(l)\geq \tl{f}(a)\geq -1 
\end{equation*}
since $\vert \tl{f}(a) \vert \leq p(a)\leq 1$ implies $-1\leq \tl{f}(a) \leq 1$, 
and condition $(2)$ holds. 

``$(2)\Ra(3)$": Assume that $x \in X\cap(A-E^+)$. Since $A=-A$, we have $x \in X\cap(-A-E^+)$, 
so there exists $-a\in -A$ and $l\in E^+$ such that $x=-a-l$, equivalently, $-x=a+l$. Thus 
$-x \in X\cap(A+E^+)$ and by condition $(2)$ we have $f(-x)\geq -1$, equivalently, 
$f(x)\leq 1$ and condition $(3)$ holds. 

``$(3)\Ra(1)$":Assume that condition $(3)$ holds. Observe that $A\subseteq (A-E^+)$ as $0\in E^+$. Direct 
checking shows that $A-E^+$ is a convex set and since $A$ is absorbing at $0$, $A-E^+$ is absorbing at $0$. Define 
\begin{equation*}
p^{'}(x):=\mathrm{inf}\{\lambda > 0 \mid x \in \lambda(A-E^+) \} 
\end{equation*}
on $E$. Then $p^{'}$ is a sublinear functional, and since $A\subseteq (A-E^+)$, $p^{'}(x) \leq p(x)$ for all $x\in E$. 

Now let $x\in X$ be an element with $p^{'}(x)\leq \lambda$ for some $\lambda > 0$, then we have $p^{'}(x/\lambda)\leq 1$. 
It follows that $x/\lambda \in X\cap (A-E^+)$ and therefore $f(x/\lambda)\leq 1$, equivalently, 
$f(x)\leq \lambda$, by the condition. Hence $f(x) \leq p^{'}(x)$ for all $x\in X$. By the Hahn-Banach theorem, 
there exists a linear extension $\tl{f}:E\ra \RR$ such that $\tl{f}(x)\leq p^{'}(x)$ for all $x\in L$. 
In particular $\tl{f}(x)\leq p(x)$ for all $x\in E$.

In order to see that $\tl{f}$ is positive, let $x\in E^+$. Then for any $\alpha >0$, $-\alpha x \in -E^+ \subseteq A-E^+$, 
so $\alpha p^{'}(-x)=p^{'}(-\alpha x)\leq 1$. Hence $-\alpha \tl{f}(x)=\tl{f}(-\alpha x) \leq p^{'}(-\alpha x)\leq 1$ and 
we obtain $-\tl{f}(x)\leq 1/\alpha$ for all $\alpha > 0$. Letting $\alpha \to \infty$, $-\tl{f}(x)\leq 0$ 
or, equivalently, $\tl{f}(x)\geq 0$ and the positivity of $\tl{f}$ is shown. 

Finally since $p$ is a seminorm we have $-\tl{f}(x)=\tl{f}(-x)\leq p(-x)=p(x)$ for all $x\in E$ and therefore 
$-p(x) \leq \tl{f}(x)$ for all $x\in E$. Together with $f(x) \leq p(x)$ for all $x\in E$, we obtain 
$\vert \tl{f}(x) \vert \leq p(x)$ for all $x\in E$ and the proof of the implication is completed.  

``$(3)\Ra(4)$": Let $x\in X$ and $y\in E$ with $x\leq y$. Then there is $l\in E_+$ 
such that $y=x+l$. Assuming $p(y)\neq 0$, we have $\displaystyle{\frac{y}{p(y)} = \frac{x}{p(y)} + \frac{l}{p(y)}}$
where $\displaystyle{\frac{x}{p(y)}  \in X}$ and $\displaystyle{\frac{l}{p(y)} \in E_+}$. Since 
$\displaystyle{p(\frac{y}{p(y)}) \leq 1}$, 
by condition $(3)$ we have $\displaystyle{f(\frac{x}{p(y)}) \leq 1}$, equivalently, $f(x) \leq p(y)$. 

If, on the other hand, $p(y) = 0$, then by similar arguments using the quotient $\displaystyle{\frac{y}{p(y)+\epsilon}}$ 
for any $\epsilon > 0$ and in the end letting $\epsilon \ra 0$, we obtain $f(x)\leq p(y)=0$.  

``$(4)\Ra(3)$": Assume that $x \in X\cap(A-E^+)$. Then there exists $l\in E_+$ such that 
$p(x+l)\leq 1$. Letting $y:=x+l$, $y\leq x$, so by condition $(4)$, $f(x)\leq p(y)\leq 1$ 
and condition $(3)$ is shown.  

``$(5)\Ra(4)$": Assume that $x\in X$ and $y\in E$ with $x\leq y$, equivalently, $-y \leq -x$. By 
condition $(5)$, $f(-(-x))\leq p(-y) = p(y)$, equivalently, $f(x)\leq p(y)$ and thus $(4)$ holds. 

``$(4)\Ra(5)$": This is similar to $(5)\Ra(4)$.

Since in particular $(1)$ holds with an extension $\tl{f}\in E^{'}_{p,+}$, clearly 
$f\in X^{'}_{p,+}$. 

Finally, assume $g\in X^{'}_{p,+}$ such that $g/\vert g \vert_{p}$ 
satisfies all the equivalent conditions $(1)$--$(5)$. By $(1)$, there exists 
a positive extension $g^{'}\colon E\ra \RR$ of $g/\vert g \vert_{p}$ such that $\vert g^{'}(x) \vert \leq p(x)$ 
for all $x\in E$. Then clearly $\tl{g}:=\vert g \vert_{p} g^{'}$ is a positive 
continuous extension of $g$ and $\vert \tl{g}(x)\vert \leq \vert g \vert_{p} p(x)$, 
so $\vert \tl{g} \vert_{p} \leq \vert g \vert_{p}$. 
Since an extension cannot decrease the norm we have $\vert \tl{g} \vert_{p} = \vert g \vert_{p}$. 
The converse direction is clear. 
\end{proof}

\begin{remark}\label{r:bnvoid}
It is interesting to observe that, Theorem \ref{t:baunami} does not require the initial 
functional $f\colon X\ra \RR$ to satisfy the natural condition $f$ is positive on $X$. Indeed, 
one can directly observe that, if $f$ is not positive on $X$, then it necessarily does 
not satisfy e.g. condition $(3)$ of the theorem: Assume that there is $x\in E^+$ and $f(x)=-\alpha < 0$. Then 
it is easy to see that $-2x/\alpha \in X\cap (A-E^+)$ with $f(-2x/\alpha)=2$, so 
the condition $(3)$ of the theorem is not satisfied. 

Similarly, the natural condition $f(x)\leq p(x)$ for all $x\in X$ is not required. Again, one can directly 
observe that, if $f$ does not have this, then it necessarily does not satisfy condition 
$(3)$, hence all conditions, of the theorem: Assume that there is $x\in X$ such that $f(x) > p(x)$. Then since 
$p$ is the Minkowski seminorm of $A$, one obtains $x\in (p(x)+\epsilon)A \subseteq (p(x)+\epsilon)(A-E^+)$ for any 
$\epsilon > 0$ and consequently $\frac{x}{p(x)+\epsilon} \in X \cap (A-E^+)$. By condition $(3)$, 
$f(\frac{x}{p(x)+\epsilon})\leq 1$, from which it follows that $f(x)\leq p(x)$, a contradiction, so such 
a functional necessarily does not satisfy the equivalent conditions in Theorem \ref{t:baunami}. 
\end{remark}

\begin{remark}\label{r:baunamihb}
Let $E$ be any real vector space. Then $E$ is a real ordered vector space considered with the 
trivial cone $E^+=\{ 0 \}$. Let $p$ be a seminorm on $E$, $X$ a subspace of $E$ 
and $f\colon X\ra \RR$ a linear functional with $f(x)\leq p(x)$ for all $x\in X$. Then conditions 
$(2)$ and $(3)$  
of the Bauer-Nachbin-Namioka Theorem are both satisfied trivially, and one recovers the existence 
of a linear extension $\tl{f}\colon E\ra \RR$ with $\vert \tl{f}(x) \vert \leq p(x)$ for all 
$x\in E$, namely, the seminorm version of the Hahn-Banach Theorem. 
\end{remark}

Given a subspace $X$ of $E$ and a seminorm $p$ on $E$, let $f\in X^{'}_{p,+}$. Then 
we call an extension $\tl{f}\in E^{'}_{p,+}$ 
with $\vert f \vert_{p} = \vert \tl{f} \vert_{p}$ a \emph{Bauer-Nachbin-Namioka extension}, 
or shortly a \emph{BNN extension}. 

The following corollary adds the equivalent condition (3) to Satz 1 from \cite{Bauer2}.
 
\begin{corollary}\label{c:bau}
Let $E$ be a real ordered vector space and $x_0\in E$ be arbitrary. Let $p$ be a seminorm 
on $E$. Then TFAE:
\begin{itemize}
\item[(1)] There exists 
a positive functional $f\colon E \ra \RR$ with $\vert f(x) \vert \leq p(x)$ for all $x\in E$ 
and $f(x_0)=p(x_0)$. 

\item[(2)] The condition $(-x_0 + A_0)\cap E^+ = \emptyset$ holds, where $A_0:=\{x\in E \mid p(x) < p(x_0) \}$ 
if $p(x_0) > 0$, and $A_0=\emptyset$ if $p(x_0)=0$. In the latter case we interpret $(-x_0 + A_0)=\emptyset$. 

\item[(3)] The condition $p(x_0)\leq p(x_0+l)$ for any $l\in E^+$ holds. 
\end{itemize}
\end{corollary}

\begin{proof}
``$(1)\Ra(2)$": If $p(x_0)=0$, then $A_0=\emptyset$, and the condition is trivially satisfied. Assume 
$p(x_0) >0$ and let $x\in A_0$. Then $f(-x_0+ x)=-f(x_0)+f(x)\leq -f(x_0)+ p(x)= -p(x_0)+ p(x) <0$. 
Since $f$ is a positive functional, $-x_0+ x \notin E^+$ and it follows that the condition is 
satisfied. \smallskip

``$(2)\Ra(3)$": If $p(x_0)=0$, then $0=p(x_0)\leq p(x_0+l)$ for any $l\in E^+$ is trivially satisfied, 
so assume $p(x_0) >0$. Let $l\in E^+$ be arbitrary and let $x=x_0+ l$, equivalently, $-x_0+ x= l$, 
so $-x_0+ x \in E^+$ and it follows that $x\notin A_0$, hence $p(x)=p(x_0+l) \geq p(x_0)$. \smallskip

``$(3)\Ra(1)$": Define the one dimensional subspace 
$X:=\{x\in E \mid x=\alpha x_0\,\,\mathrm{for\,\,some}\,\,\alpha\in\RR  \}$ of $E$ and a linear 
functional $f\colon X\ra \RR$ on $X$ by $f(x):=\alpha p(x_0)$; note that, $f(x_0)=p(x_0)$. Let $A:=\{ x\in E \mid p(x)\leq 1 \}$. 
Let $\alpha x_0 = x\in X\cap(A-E^+)$ be arbitrary. In order to use Theorem \ref{t:baunami} we look at two cases: \smallskip

\emph{Case I:} If $\alpha \leq 0$, then $f(x)=\alpha p(x_0) \leq 0 < 1$. \smallskip

\emph{Case II:} If $\alpha >0$, we claim that, $f(x)\leq 1$, equivalently, $p(x_0) \leq 1/\alpha$ must hold: 
Since $\alpha x_0 \in A-E^+$, there exists $l\in E^+$ such that $\alpha x_0 + l \in A$, so 
$p(\alpha x_0 + l) \leq 1$, equivalently, $p(x_0+ (1/\alpha)l) \leq 1/\alpha$. Since 
$(1/\alpha)l \in E^+$, it follows that $p(x_0)\leq 1/\alpha$. \smallskip

By Theorem \ref{t:baunami} 
we conclude that a positive extension of $f$, again denoted by $f$, with $\vert f(x) \vert \leq p(x)$ for all 
$x\in E$ exists, so (1) is shown. \smallskip
\end{proof}

Given an ordered space $E$ with a seminorm $p$, we call a positive functional 
$f\in E^{'}_{p,+}$ a \emph{$p$-state} or a \emph{state} if $p$ is clear, if $\vert f \vert_{p}=1$. 
Note that if $f$ is a $p$-state, then in particular we have $\vert f(x) \vert\leq p(x)$ for all $x\in E$. 
Corollary \ref{c:bau} gives equivalent conditions for the existence of a certain type 
of state. In the rest of the article we will see classes of ordered spaces with a 
locally convex topology where the set of states is nonempty. 

\begin{remarks}
\begin{itemize}
\item[(1)]In Corollary \ref{c:bau} it is clear that if $p(x_0)=0$, then all the equivalent conditions are trivially 
satisfied and in particular a positive linear extension in (1) is the zero functional. 

\item[(2)]As in Remark \ref{r:baunamihb} let $E$ be a real vector space with a seminorm $p$, 
considered as a real ordered vector space with $E^+:=\{ 0 \}$. 
In this case (3) of Corollary \ref{c:bau} holds trivially for any $x_0\in E$, 
and one recovers the existence of a linear functional $f\colon E\ra \RR$ with 
$f(x_0)=p(x_0)$ and $f(x)\leq p(x)$ for all $x\in E$, a well known corollary of the Hahn-Banach Theorem.  
\end{itemize} 
\end{remarks}

By considering $-x_0$ in place of $x_0$, we obtain the following corollary of Corollary \ref{c:bau}. 
Alternatively, one could use condition $(2)$ of the Bauer-Namioka theorem with similar ideas 
as in the proof of Corollary \ref{c:bau}.

\begin{corollary}\label{c:bau2}
Let $E$ be a real ordered vector space and $x_0\in L$ be arbitrary. Let $p$ be a seminorm 
on $E$. Then TFAE:
\begin{itemize}
\item[(1)] There exists 
a positive functional $f\colon E \ra \RR$ with $\vert f(x) \vert \leq p(x)$ for all $x\in E$ 
and $f(x_0)=-p(x_0)$. 

\item[(2)] The condition $(x_0 + A_0)\cap E^+ = \emptyset$ holds, where $A_0:=\{x\in E \mid p(x) < p(x_0) \}$ 
if $p(x_0) > 0$, and $A_0=\emptyset$ if $p(x_0)=0$. In the latter case we interpret $(x_0 + A_0)=\emptyset$. 

\item[(3)] The condition $p(x_0)\leq p(x_0-l)$ for any $l\in E^+$ holds. 
\end{itemize}
\end{corollary}

By combining the statements of Corollaries \ref{c:bau} and \ref{c:bau2}, more precisely, 
using the equivalence of items $(1)$ and $(3)$ of them, we obtain the following corollary.

\begin{corollary}\label{c:bau3}
Let $E$ be an ordered vector space and let $p$ be a seminorm 
on $E$. Assume that all elements $x\in E$ satisfy at least one of the following: 
$p(x+l)\geq p(x)$ for all $l\in E_+$ or $p(x-l)\geq p(x)$ for all $l\in E_+$. Then 
for any $x_0\in E$ there exists a positive functional $f\colon E \ra \RR$ 
with $\vert f(x) \vert \leq p(x)$ for all $x\in L$ 
and $\vert f(x_0) \vert=p(x_0)$.  
\end{corollary}

In the case the seminorm $p$ is increasing, the following corollary of Corollary \ref{c:bau} can be obtained.

\begin{corollary}\label{c:bauay}
Let $E$ be a real ordered vector space and assume that $p$ be an increasing seminorm on $E$. 
Then for any $x_0\in E^+$ there exists a positive 
functional $f\colon E \ra \RR$ with $f(x_0)=p(x_0)$ and $\vert f(x) \vert \leq p(x)$ for all $x\in E$.  
\end{corollary}

\begin{proof}
Since $x_0\in E^+$ and $p$ is increasing, $p(x_0)\leq p(x_0+ l)$ for all $l\in E^+$ 
and condition (3) of Corollary \ref{c:bau} is satisfied, so there exists a positive 
functional with the asserted properties by Corollary \ref{c:bau}. 
\end{proof}

We record the following corollary of Corollary \ref{c:bauay} for future reference.

\begin{corollary}\label{c:bauay2}
Let $E$ be an ordered vector space and assume that $p$ is an increasing seminorm on $E$. 
Then for any $x_0\in E^+$ we have 
\begin{equation*}
p(x_0)=\underset{0\leq f\leq p}{\mathrm{sup}}\,\vert f(x_0) \vert= \underset{0\leq f\leq p}{\mathrm{sup}}\,f(x_0)
\end{equation*}
where the supremum is taken over all $f\in E_{+}^{'}$ with $\vert f(x) \vert\leq p(x)$ for all $x\in E$. 
\end{corollary}

\section{Schaefer's Representation Theorem} \label{s:isoschrep}

In this section we provide a version of Schaefer's representation theorem, in which the representation is isometric on the positive cone. 
Schaefer's representation theorem is originally due Schaefer in \cite{Schaefer2}, (5.1), where the closedness condition of the cone was missing, see also \cite{Schaefer1}, 4.4 and \cite{Riedl}, Theorem 4.1 for correct versions. On the other hand 
we could not find the version that we state in this article elsewhere, though we believe it must be folklore. Lemma 2.5. 
of \cite{Davies} for an ordered Banach space with a certain regularity condition on its dual is the closest statement we could find; 
see also \cite{Wittstock} for such spaces.  

See also the recent article \cite{deBruyn1} where in Theorem B characterizations of injective, positive and continuous representability of a real ordered topological vector space $E$ on the space $C(X)$ where $X$ is locally compact are obtained.

Throughout the rest of the article we will assume that $E$ is an ordered space with a locally convex topology, with the cone $E^+$. Notice that 
we do not assume closedness or normality of the cone $E^+$. We will consider 
such a space with a specified \emph{calibration}, that is, a specified family $\{ p_{\alpha} \}_{\alpha\in A}$ of 
continuous seminorms generating the topology of $E$. We always 
assume that the family $\{ p_{\alpha} \}_{\alpha\in A}$ Hausdorff separates $E$, 
that is, $p_{\alpha}(x)=0$ for all $\alpha\in A$ implies $x=0$. Throughout the remainder of the article we will refer to such a 
space $(E, E_+, \{ p_{\alpha} \}_{\alpha\in A})$ as a \emph{calibrated ordered space}. The family 
$\{ p_{\alpha} \}_{\alpha\in A}$ is called \emph{saturated}, if it is 
closed under the pointwise maximum seminorm, 
i.e. for any $\{p_{\alpha_i}\}_{i=1}^n$  from the family, 
the seminorm defined by $p(x):=\mathrm{max}_{i=1}^n p_{\alpha_i}(x)$ for any $x\in X$ 
is also in the family. Note that if a given family $\{ p_{\alpha} \}_{\alpha\in A}$ 
does not have this property, we can always saturate it by adding pointwise maximum seminorms 
for all finite families and clearly this operation does not disturb the topology. 
Moreover, this operation preserves some other properties of the seminorms in the family in different 
settings, as we will see. 

We denote the topological dual space of $E$ by $E^{'}$. Let $f\colon E\ra \RR$ be a linear functional in $E^{'}$ and 
note that if the generating family $\{ p_{\alpha} \}_{\alpha\in A}$ is saturated, 
there exists $\alpha\in A$ and a constant $C\geq 0$ such that $\vert f(x) \vert \leq Cp_{\alpha}(x)$ for all $x\in E$. 
For a fixed $\alpha\in A$, we denote the set of all such $f$ by $E^{'}_{p_{\alpha}}$ or more 
concisely by $E^{'}_{\alpha}$, similar to the notation used in Section \ref{s:nps}. It is 
clear that $E^{'}_{\alpha}$ is a vector subspace of $E^{'}$. 
If $f\in E^{'}_{\alpha}$, let $\vert f \vert_{\alpha}$ be the infimum of all $C\geq 0$ 
with $\vert f(x) \vert \leq Cp_{\alpha}(x)$ for all $x\in E$. It is easy to see that $\vert \cdot \vert_{\alpha}$ 
is a norm on $E^{'}_{\alpha}$. 
By $E^{'}_{\alpha,+}$ we denote the subset of $E^{'}_{\alpha}$ containing positive functionals. 
Finally, by $E^{'}_{\alpha,1}$ we denote the set of all $f$ with $\vert f \vert_{\alpha} \leq 1$ 
and by $E^{'}_{\alpha,+,1}$ the subset of $E^{'}_{\alpha,1}$ containing positive functionals. 

An important example of a calibrated ordered space is $C(X)$: Let $X$ be a Hausdorff topological space, and let $C(X)$ be the space 
of all continuous real valued functions on $X$. Then $C(X)$ is an ordered vector space with a natural cone consisting 
of everywhere positive or zero functions. It is easy to see that $C(X)$ is an Archimedean 
ordered space and hence also almost Archimedean. Given a collection of compact subsets $\{ K_{\alpha} \}_{\alpha\in A}$, 
we take $p_{\alpha}$ to be the seminorm of supremum over $K_{\alpha}$. 
If the specified compact sets are closed under finite unions then the corresponding seminorm family $\{ p_{\alpha} \}_{\alpha\in A}$ 
is closed under finite pointwise maximums. Notice that, since a finite union 
of compact subsets is compact, we can always saturate 
the collection $\{ K_{\alpha} \}_{\alpha\in A}$ by 
adding all finite unions without altering the topology. 
On the other hand, if $X=\bigcup_{\alpha\in A} K_{\alpha}$, it is then easy to see that 
the locally convex topology is Hausdorff, 
that is, if $p_{\alpha}(f)=0$ for all $\alpha\in A$, $f=0$ in $C(X)$. 

We will denote $C(X)$ with the indicated structure by $(C(X),C(X)_+,\{ K_{\alpha}\}_{\alpha\in A} )$ or 
shortly by $C(X)$, if the rest is clear from the text. In the next section we will see that 
a certain class of calibrated ordered spaces always arise in this way. 

We recall that, the \emph{compact open} topology on $C(X)$ is the topology of uniform convergence on compact subsets of $X$, 
in other words, a net in $C(X)$ converges if and only if it converges uniformly on compact subsets of $X$. 

We state the following fact, which is a well known consequence of the Hahn - Banach Theorem, see e.g. \cite{Conway}, Theorem 3.7, 
as well as \cite{Eidelheit}.  

\begin{lemma}\label{lhb}
Let $X$ be a real topological vector space and $A$ and $B$ be two nonempty disjoint convex subsets of $X$ with $A$ open. Then 
there exists a continuous linear functional $f\colon X\to \RR$ and a number $\omega\in\RR$ such that $f(a)<\omega$ for all 
$a\in A$ and $f(b)\geq \omega$ for all $b\in B$.  
\end{lemma}

The following type of positive Hahn-Banach type separation lemma is well known, see e.g. Theorem 2.13(4) in \cite{Aliprantis}, 
Lemma 1.2 in \cite{Kelley}, as well as Corollary 4.2 in \cite{Namioka}. We provide a proof. 

\begin{theorem}\label{t:rphbs}
Let $(E, E_+, \{ p_{\alpha} \}_{\alpha\in A})$ be a calibrated ordered space where 
the cone $E_+$ is closed. If a point $x_0\in E$ satisfies $f(x_0)\geq 0$ for all $f\in E^{'}_{+}$, then $x_0\geq 0$.
\end{theorem}

\begin{proof}
Assume $x_0\notin E_+$. Then since $E_+$ is closed, there exists an open and convex set $A\subseteq E$ containing $x_0$ such that $A\cap E_+ = \emptyset$. 
Since $E_+$ is convex, by Lemma 
\ref{lhb} there exists a real linear continuous functional $f\colon E\to \RR$ and $\omega\in\RR$ such that $f(y)< \omega$ for all $y \in A$ and 
$f(z)\geq \omega$ for all $z\in E_+$. As $0\in E_+$, $0=f(0)\geq \omega$ so we have $f(x_0)< 0$. Now let $x_1\in E_+$ be arbitrary. 
Then $nx_1 \in E_+$ for all $n\in \NN$, so 
$f(nx_1)=n f(x_1)\geq \omega$. Dividing by $n$ and letting $n\to\infty$, we get $f(x_1)\geq 0$. Hence $f(E_+)\subseteq \RR_+$ and 
$f\in E^{'}_{+}$, a contradiction. Thus $x_0\in E_+$. 
\end{proof}

The following fact is well known, first stated and proved in \cite{Krein2} for the normed case. 
See also \cite{Schaefer1}, p.218 and p.219, \cite{Aliprantis}, Theorem 2.26 and \cite{Riedl}, 
Corollary 3.2 and Corollary 3.3. We give 
a statement which refers to calibrated ordered spaces. 

\begin{theorem}\label{t:krein}
Let $(E, E_+, \{ p_{\alpha} \}_{\alpha\in A})$ be a calibrated ordered space with 
the family $\{ p_{\alpha} \}_{\alpha\in A})$ consisting of increasing seminorms, equivalently, 
the cone $E_+$ is normal, see Section \ref{s:nps}. Let $f\in E^{'}$. Then there exists $u,v\in E^{'}_{+}$ 
such that $f=u-v$. 
\end{theorem}

Given a finite family $\{p_{\alpha_i}\}_{i=1}^n$ of increasing seminorms on an 
ordered space $E$, it is easy to check that the pointwise maximum seminorm is also increasing. 
Therefore it makes sense to consider a calibrated ordered space 
$(E, E_+, \{ p_{\alpha} \}_{\alpha\in A})$ with a saturated increasing family of seminorms. 

The following theorem is a version of Schaefer's representation theorem, with the clause 
that the representation can be done with a map which is ``isometric" on the set of all positive elements. 

\begin{theorem}\label{t:schaf}
Let $(E, E_+, \{ p_{\alpha} \}_{\alpha\in A})$ be a calibrated ordered space where 
the family $\{ p_{\alpha} \}_{\alpha\in A}$ is closed under pointwise maximum. Then TFAE:
\begin{itemize}
\item[(1)] The cone $E^+$ of $E$ is closed and the family $\{ p_{\alpha} \}_{\alpha\in A}$ 
consists of increasing seminorms. 

\item[(2)] There exists a locally compact Hausdorff space $X$ and an ordered subspace $W$ 
of $C(X)$ such that there is a linear bipositive bijection $\phi\colon E\ra W$, 
which is continuous and its inverse $\phi^{-1}\colon W \ra E$ is also continuous. In addition, 
the space $X$ and the map $\phi$ can be constructed in such a way that, there exists a family 
of pairwise disjoint compact subsets $\{K_{\alpha}\}_{\alpha\in A}$ of $X$, defining 
the topology of $C(X)$ by the supremum on $K_{\alpha}$ seminorms, with the properties 
$X=\bigcup_{\alpha\in A} K_{\alpha}$ and any compact subset $K\subseteq X$ 
is contained in a finite union of the specified compact subsets $\{ K_{\alpha} \}_{\alpha\in A}$. Moreover 
we have $\vert \phi(x) \vert_{K_{\alpha}}= p_{\alpha}(x)$ for all $x\in E^+$. 
\end{itemize}
\end{theorem}

\begin{proof}
``$(1)\Ra(2)$": The first and second assertions are in the content of Schaefer's representation theorem, see e.g. the proof of 4.4 in \cite{Schaefer1}. 
However, since we have to refer to a slightly more specific construction of $X$ for the proof of the remaining assertion, 
we provide full details. For any $p_{\alpha}$, with $\alpha\in A$ let 
$B_{\alpha}:=E_{\alpha,+,1}=\{ f\in E_{+}^{'} \mid \vert f(x) \vert \leq p_{\alpha}(x)\,\,\mathrm{for\,\,all\,\,}x\in E\}$. By 
the Alaoglu-Bourbaki theorem, each $B_{\alpha}$ is compact under the weak-$*$ topology. Define a set $X$ 
to be the disjoint union of all $B_{\alpha}$, e.g. let 
$X:=\cup_{\alpha\in A} B_{\alpha}\times \{ \alpha \}$. We consider $X$ with the disjoint union topology, 
i.e. by definition a subset is open if and only if its intersection with all $B_{\alpha}$ are
open. Note that, a base of neighbourhoods of any $(f,\alpha)$ is given by sets of type 
$\{ (g,\beta) \mid \vert f(x)-g(x) \vert < \epsilon \}$ for all $x\in E$ and $\epsilon > 0$. In particular, a net $(f_j,\alpha_j)_{j\in J}$ 
converges to an element $(f,\alpha)$ if and only if $\displaystyle{\vert f_j(x) - f(x) \vert \underset{j\in J}{\ra} 0}$ for all $x\in E$. 

We observe that $X$ is a locally compact space. Indeed, for any $(f,\alpha)\in X$ we have $f\in B_{\alpha}$, 
and $(f,\alpha)\in B_{\alpha}\times \{\alpha \}$, where $B_{\alpha}$ 
is compact. So $B_{\alpha}\times \{\alpha \}$ is compact and open by construction, hence it is a compact neighbourhood of $(f,\alpha)$. 
Also, $X$ is Hausdorff. If $(f, \alpha) \neq (g,\beta)$ with $\alpha\neq \beta$, then consider disjoint neighbourhoods 
$B_{\alpha} \times \{\alpha \}$ and $B_{\beta} \times \{\beta \}$. On the other hand, if $\alpha=\beta$ and $f\neq g$, then 
there is $x\in E$ such that $f(x)\neq g(x)$. By definition of the weak-$*$ topology, one can find disjoint neighbourhoods 
of $(f,\alpha)$ and $(g,\beta)$. 

Now let us define a map $\phi\colon E\to C(X)$ by $\phi(x)(f,\alpha):=f(x)$. If $(f_j,\alpha_j)_{j\in J}$ is a net 
converging to $(f,\alpha)$, then $\phi(x)(f_j,\alpha_j)=f_j(x)\underset{j\in J}{\ra} f(x)=\phi(x)(f,\alpha)$ for any $x\in E$, 
so $\phi$ is indeed valued in $C(X)$. It is easy to check that $\phi$ is linear. In order to see that $\phi$ is an injection, 
let us show that $\mathrm{ker}(\phi)=\{0\}$. For a contradiction, assume that there exists 
$x\neq 0$ in $E$ such that $0=\phi(x)(f,\alpha)=f(x)$ for all $(f,\alpha)\in X$. Since $E$ is Hausdorff separated 
there exists a functional $s\in E^{'}$ such that 
$s(x)\neq 0$. By Theorem \ref{t:krein} there exist positive functionals $u,v\in E_{+}^{'}$ with $s=u-v$. Since 
$\{ p_{\alpha} \}_{\alpha\in A}$ is saturated there exists $p_{\beta}, p_{\gamma}$ and $c_1, c_2 > 0$ 
such that $\vert u(x) \vert\leq c_1 p_{\beta}(x)$ 
and $\vert v(x) \vert\leq c_2 p_{\gamma}(x)$ for all $x\in E$, so that $(u/c_1,\beta), (v/c_2,\gamma)\in X$. It follows that 
$u(x)=0$ and $v(x)=0$, a contradiction. Thus $\phi$ is an injection. 

We now show that $\phi$ is a bipositive map. If $x\in E_{+}$, then $\phi(x)(f,\alpha)=f(x)\geq 0$ for all $(f,\alpha)\in X$, 
so $\phi(x)\geq 0$ in $C(X)$ by definition of positivity in $C(X)$. Thus $\phi$ is a positive map. If, on the other hand, 
$\phi(x)\geq 0$ for some $x\in E$, by definition of $\phi$ we obtain $\phi(x)(f,\alpha)=f(x)\geq 0$ for all $(f,\alpha)\in X$ 
and thus $f(x)\geq 0$ for all $f\in E^{'}_{+}$. By Theorem \ref{t:rphbs} we have $x\geq 0$ and we conclude that $\phi$ 
is a bipositive map. 

In order to see that $\phi$ is continuous, let $K\subset X$ be a compact set. Then there exists finitely many 
$(B_{\alpha_j}\times\{\alpha_j\})_{j=1}^n$ covering $K$. Then we have
\begin{align*}
\vert \phi(x) \vert_{K}&=\underset{(f,\alpha)\in K}{\mathrm{sup}}\,\vert \phi(x)(f,\alpha) \vert \\
&=\underset{(f,\alpha)\in K}{\mathrm{sup}}\,\vert f(x) \vert \\
&\leq \underset{j=1,\cdots,n}{\mathrm{max}}\,p_{\alpha_j}(x)
\end{align*}
for all $x\in E$ and continuity of $\phi$ is shown. 

We show that the inverse mapping $\phi^{-1}\colon \phi(E)(\subseteq C(X))\ra E$ is continuous.  
If $\phi(x_j)\ra\phi(x)$, where $(x_j)_{j\in J}$ is a net in $E$ and $x\in E$, is a 
convergent net in $\phi(E)$, $(\phi(x_j))_{j\in J}$ is uniformly convergent to $\phi(x)$ on compact subsets of 
$X$ by definition of the topology on $C(X)$.  In particular, $\phi(x_j)\ra\phi(x)$ uniformly on each 
compact set $B_{\alpha}\times \{ \alpha \}$, equivalently, $f(x_j)\ra f(x)$ uniformly on all 
$f\in B_{\alpha}$. Since $\{(p_{\alpha})_{\alpha\in A}\}$ is generating, it follows that $f(x_j)\ra f(x)$ uniformly on 
equicontinuous subsets of $E_{+}^{'}$. Now we use the fact that, convergence of a net in $E$ is characterized by the uniform convergence 
of it on 
equicontinuous subsets of $E_{+}^{'}$, see 3.3, 
Corollary 1 in \cite{Schaefer1}, to conclude that $x_j\ra x$ in $E$. Note that, reversing the steps, this argument 
would show the continuity of $\phi$ as well. 

This concludes the proof of the first assertion. To prove the remaining assertions, let $K_{\alpha}:=B_{\alpha}\times\{ \alpha \}$ 
for any $\alpha\in A$. Then each $K_{\alpha}$ is compact, the sets $\{ K_{\alpha} \}_{\alpha\in A}$ are pairwise disjoint 
and $X=\bigcup_{\alpha\in A} K_{\alpha}$. Now given an arbitrary $\alpha\in A$ and $x\in E^{+}$ we have 
\begin{align*}
\vert \phi(x) \vert_{K_{\alpha}}&= \underset{(f,\alpha)\in K_{\alpha}}{\mathrm{sup}}\,\vert \phi(x)(f,\alpha) \vert \\
&=\underset{(f,\alpha)\in K_{\alpha}}{\mathrm{sup}}\,\vert f(x) \vert \\
&=\underset{0\leq f\leq p_{\alpha}}{\mathrm{sup}}\,\vert f(x) \vert \\
&=p_{\alpha}(x)
\end{align*} 

where the last equality follows by Corollary \ref{c:bauay2} and the theorem is proven. 

``$(2)\Ra(1)$":By the first and second assertions, this implication is included in the usual form of Schaefer's theorem, see e.g. 
the proof of 4.4 in \cite{Schaefer1}.
 
\end{proof} 

If in Theorem \ref{t:schaf} we insist on having a calibrated ordered space 
$(E, E_+, \{ p_{\alpha} \}_{\alpha\in A})$ where the specified family of seminorms 
is not necessarily saturated under the pointwise maximum, then the following version 
is obtained. 

\begin{corollary}\label{c:schaf}
Let $(E, E_+, \{ p_{\alpha} \}_{\alpha\in A})$ be a calibrated ordered space. Then TFAE:
\begin{itemize}
\item[(1)] The cone $E^+$ of $E$ is closed and the family $\{ p_{\alpha} \}_{\alpha\in A}$ 
consists of increasing seminorms. 

\item[(2)] There exists a locally compact Hausdorff space $X$ and an ordered subspace $W$ 
of $C(X)$ such that there is a linear bipositive bijection $\phi\colon E\ra W$, 
which is continuous and its inverse $\phi^{-1}\colon W \ra E$ is also continuous. Moreover, 
the space $X$ and the map $\phi$ can be constructed in such a way that, there exists a family 
of pairwise disjoint compact subsets $\{K_{\beta}\}_{\beta\in B}$ of $X$, 
where $A\subseteq B$, 
the topology of $C(X)$ is defined by the supremum on $\{K_{\alpha}\}_{\alpha\in A}$ seminorms, with the property 
$\vert \phi(x) \vert_{K_{\alpha}}= p_{\alpha}(x)$ for all $x\in E^+$ and $\alpha\in A$,  
$X=\bigcup_{\beta\in B} K_{\beta}$ and any compact subset $K\subseteq X$ 
is contained in a finite union of the compact subsets $\{ K_{\beta} \}_{\beta\in B}$. 
\end{itemize}
\end{corollary}

\begin{proof}
Consider the calibrated ordered space $(E, E_+, \{ p_{\beta} \}_{\beta\in B})$ obtained by 
saturating the family $\{ p_{\alpha} \}_{\alpha\in A}$ by all pointwise maximum seminorms; 
in particular, $A\subseteq B$. Then applying Theorem \ref{t:schaf}, we obtain 
a locally compact Hausdorff space $X$, an ordered subspace $W$ of $C(X)$ and 
a map $\phi\colon E\ra W$ with the stated properties, since the topology of $E$ 
with the saturated family remains the same. 
\end{proof}

Specializing to a normed ordered space, the following corollary is obtained. 

\begin{corollary}\label{c:nschaf}
Let $E$ be a real ordered normed space. Then TFAE:
\begin{itemize}
\item[(1)] The cone $E^+$ of $L$ is closed and normal. 

\item[(2)] There exists a compact Hausdorff space $K$ such that there is a linear bipositive bijection $\phi\colon E\ra W$, 
which is continuous and its inverse $\phi^{-1}\colon W \ra E$ is also continuous, 
where $W$ is an ordered subspace of $C(K)$. Moreover, 
we have $\vert \phi(x) \vert_{K}= \Vert x\Vert_{E}$ for all $x\in E^+$. 
\end{itemize}
\end{corollary} 

\section{Characterizations of Isometric Representability of a Calibrated Ordered Space}\label{s:cirpolcs}

In this section we will state and prove theorems \ref{t:isorepcx}, \ref{t:isorepsnfl}, \ref{t:isorepstr} and \ref{t:isorepquot}, 
the main theorems of the article, 
all of which list equivalent conditions for the isometrically isomorphic representations 
of calibrated ordered spaces on a space of continuous functions in the sense below. Some of the items 
listed are well known in the normed case, and some of them appear to be new even in the normed case, 
see Remark \ref{r:refofisorep} below. 

Two calibrated ordered spaces $(E, E_+, \{ p_{\alpha} \}_{\alpha\in A})$ and 
$(F, F_+, \{ q_{\beta} \}_{\beta\in B})$ are \emph{isometrically bipositively isomorphic}, 
or shortly \emph{isometrically isomorphic} if $A=B$ and there exists a linear bijection 
$\phi\colon E\ra F$ such that $\phi$ is bipositive, that is, 
$x \in E_+ \iff  \phi(x)\in F_+$ and for any $\alpha \in A$, $q_{\alpha}(\phi(x))=p_{\alpha}(x)$ 
for all $x\in E$. Notice that, the isometry condition on each $\alpha\in A$ guarantees 
continuity of both $\phi$ and its inverse $\phi^{-1}$. 

Following \cite{Nachbin}, given a calibrated ordered space $(E, E_+, \{ p_{\alpha} \}_{\alpha\in A})$ 
we say that an element $x\in E$ is \emph{semi positive} if $p_{\alpha}(x+l)\geq p_{\alpha}(x)$ for any $\alpha\in A$ 
and $l\in E_+$. Similarly, $x\in E$ is called \emph{semi negative} if 
$p_{\alpha}(x-l)\geq p_{\alpha}(x)$ for any $\alpha\in A$ and $l\in E_+$. 

Given a topological space $X$, with notation as in Section \ref{s:isoschrep}, recall that the 
compact open topology on $C(X)$ is the topology of uniform convergence on all compact subsets. In general, 
the topology on $C(X)$ as a calibrated ordered space may be different than the compact open topology on $C(X)$. 
\begin{example}
Consider $X=\RR$ with a collection of compact subsets the collection of all finite subsets. Then the topology 
obtained on $C(\RR)$ is the topology of pointwise convergence, which is different then the compact open topology, equivalently, 
the topology of uniform convergence on compact subsets of $\RR$.
\end{example}

In some cases, however,  
the topology of a calibrated ordered space and the compact open topology coincide, e.g. if the space 
$X$ is locally compact Hausdorff and one considers the collection of all compact subsets of $X$. 
Since any compact subset $K$ is contained in a finite union of the fundamental subsets $B_{\alpha}\times\{\alpha\}$, for the 
locally compact Hausdorff space $X$ constructed in the proof of Theorem \ref{t:schaf} the compact open topology 
and the topology obtained by the collection $\{B_{\alpha}\times \{\alpha\}  \}_{\alpha\in A}$ coincide as well. Also, 
in the special case when the space $X$ is compact and Hausdorff, we consider $C(X)$ with its usual supremum on $X$ 
norm topology, which clearly is the compact open topology on $C(X)$. See also Remark 2.11 in \cite{Asadi}. 

Note that since 
a locally compact Hausdorff space $X$ is a \emph{k-space}, i.e. a subset is open if and only if 
its intersection with every compact subset is open in the subset topology, and since a locally compact Hausdorff space in addition  
is completely regular, by e.g. \cite{Frago} p.104, $C(X)$ with 
the compact-open topology is complete; see
also \cite{Collins} and \cite{Michael}. Therefore, if we assume that $X$ has the property that, 
any compact subset is contained in a finite union of the specified compact subsets $\{ K_{\alpha} \}_{\alpha\in A}$, 
then as in the above example, the locally convex topology generated by $\{ K_{\alpha} \}_{\alpha\in A}$ is 
the compact-open topology, and in that case, $C(X)$ is a complete locally convex space. 
Moreover, by e.g. 
\cite{Frago} p.104, $\cC(X):=C(X)+iC(X)$, the space of all 
$\CC$-valued continuous functions on $X$, with topology given by the supremum 
on $K_{\alpha}$ seminorms is a commutative locally $C^*$-algebra. Clearly, $\cC(X)^{h}=C(X)$. 

\begin{remark}\label{r:indfl}
Assume that two calibrated ordered spaces $(E, E_+, \{ p_{\alpha} \}_{\alpha\in A})$ and 
$(F, F_+, \{ q_{\beta} \}_{\beta\in B})$ are isometrically isomorphic via $\phi\colon E\ra F$ and let $f\in E^{'}$. 
Then it is easy to see that the functional $g(y):=f(\phi^{-1}(y))$ for each $y\in F$ is in $F^{'}$. In fact, 
$\vert f \vert_{\alpha}=\vert g \vert_{\alpha}$ if $f\in E^{'}_{\alpha}$ for some $\alpha\in A$. If furthermore 
$f\in E^{'}_{\alpha,+}$, then in addition $g\in E^{'}_{\alpha,+}$ with the same $\alpha$-norm. 
\end{remark}

Given a locally ordered space $(E,\{ C_{\alpha} \}_{\alpha\in A})$ with a local order unit $e$, see Section \ref{s:nps}, 
for any $\alpha\in A$ and $x\in E$ define 
$p_{\alpha}(x):=\mathrm{inf} \{\lambda >0 \mid -\lambda e \leq_{\alpha} x \leq_{\alpha} \lambda e  \} $. 
It is easy to see that $\{p_{\alpha}\}_{\alpha\in A}$ is a family of seminorms. We will always consider 
a locally ordered space $(E,\{ C_{\alpha} \}_{\alpha\in A})$ with a local order unit $e$ with 
this family of seminorms. If the space is Archimedean, then the cone $E_+$ is automatically closed: 
Let $(v_i)_{i\in I}\subseteq E_+$ be a net with $v_i\ra v$. Then for any $\alpha\in A$ we have $p_{\alpha}(v_i-v)\ra 0$. 
The definition of the seminorm implies that, for a fixed $\alpha\in A$ and any $\lambda > 0$ 
there exists $i\in I$ such that $\lambda e \geq_{\alpha} v_i - v$, equivalently, $\lambda e + v \geq_{\alpha} v_i$. 
Thus $\lambda e + v \in C_{\alpha}$ for all $\lambda > 0$ and by the Archimedean property $v\in C_{\alpha}$. 
Since this can be done for all $\alpha \in A$, $v\in \cap_{\alpha\in A} C_{\alpha} = E_+$ and $E_+$ 
is closed.   

\begin{example}\label{r:cxlocou}
Let $X$ be a locally compact Hausdorff space. 
Consider the calibrated ordered space $(C(X),C(X)_+,\{ K_{\alpha} \}_{\alpha\in A})$ with a specified family of compact subsets 
$\{ K_{\alpha} \}_{\alpha\in A}$ such that $\cup_{\alpha\in A} K_{\alpha} = X$. Taking corresponding wedges 
$C_{\alpha}:=\{f\in C(X) \mid f_{|K_{\alpha}}\geq 0 \}$ and as an order unit the constant function $1$, 
it is easy to see that $C(X)$ an L-Archimedean local order unit space. 
Notice that the wedges $\{ C_{\alpha} \}_{\alpha\in A}$ does not satisfy the downward filtered 
property, see Remark \ref{r:dosidef}. 
\end{example}

A calibrated ordered space $(E, E_+, \{ p_{\alpha} \}_{\alpha\in A})$ is called  \emph{locally Riesz space} 
if the pair $(E, E_+)$ is a Riesz space; i.e. for any pair $x,y\in E$ their supremum $x\vee y$ exists in $E$, 
and for each $\alpha\in A$ $p_{\alpha}$ is a Riesz seminorm, that is, $p_{\alpha}(\vert x \vert)=p_{\alpha}(x)$ 
for any $x\in E$. Here as usual $\vert x \vert= x \vee (-x)$. See p.14 of \cite{Aliprantis}. 

If the seminorm $p_{\alpha}$ has the further property that, $p_{\alpha}(x \vee y) = p_{\alpha}(x) \vee p_{\alpha}(y)$ 
for any pair $x,y\in E_+$, then $p_{\alpha}$ is called an M-seminorm. See \cite{Schaefer1}, also \cite{vanGaans}. 
Clearly, an M-seminorm is increasing.  

Given two M-seminorms $p$ and $q$ on a Riesz space $(E, E_+)$, let $r$ be the pointwise maximum seminorm. 
Then $r(\vert x \vert)=\mathrm{max}\{p(\vert x \vert), q(\vert x \vert)\}= \mathrm{max}\{p(x), q(x)\}=r(x)$ 
for any $x\in E$, so $r$ is a Riesz seminorm. Also,
\begin{align*}  
r(x\vee y)&=\mathrm{max}\{p(x\vee y), q(x\vee y) \} \\
&= \mathrm{max}\{p(x)\vee p(y)), q(x)\vee q(y)) \} \\
&= \mathrm{max}\{p(x)\vee q(x))\}\vee \mathrm{max}\{p(y)\vee q(y)) \} \\
&=r(x) \vee r(y)
\end{align*}
for any $x,y\geq 0$, so $r$ is an M-seminorm as well. Consequently, a locally Riesz space 
$(E, E_+, \{ p_{\alpha} \}_{\alpha\in A})$ where the specified seminorm collection consists of M-seminorms 
can always be saturated with pointwise maximum seminorms without altering the other structures.   
 
\begin{example}\label{r:cxriesz}
If we take $X$ a locally compact Hausdorff space, 
then the calibrated ordered space $(C(X),C(X)_+,\{ K_{\alpha} \}_{\alpha\in A})$ with a specified family of compact subsets 
$\{ K_{\alpha} \}_{\alpha\in A}$, with $x\vee y$ taken to be the pointwise maximum function of $x$ and $y$ 
as usual and supremum seminorms on the compact sets $K_{\alpha}$, 
is a locally Riesz space. Clearly the seminorms are M-seminorms. As we will see in Theorem \ref{t:isorepstr}, 
this is a prototype of locally Riesz spaces with topology given by M-seminorms. 
\end{example}

We now state and prove the following extension theorem, which is an application of Schaefer's construction in 
Section \ref{s:isoschrep} and the Bauer-Nachbin-Namioka Theorem.  

\begin{theorem}\label{t:posext}
Let $(E, E_+, \{ p_{\alpha} \}_{\alpha\in A})$ be a calibrated ordered space with a closed cone $E_+$ . Assume 
that $E$ has the property that, 
for any $\alpha\in A$ and  
$x_0\in E$ there exists $f\in E_{\alpha,+}^{'}$ with $\vert f(x_0) \vert=p_{\alpha}(x_0)$ and $\vert f(x) \vert \leq p_{\alpha}(x)$ for all $x\in E$. 
Then there exists a locally compact Hausdorff space $X$, 
a calibrated ordered space $(C(X),C(X)_+,\{ K_{\alpha} \}_{\alpha\in A})$ 
such that $\cup_{\alpha\in A}K_{\alpha}=X$ and any compact subset $K$ 
is contained in a finite union from the family $\{ K_{\alpha} \}_{\alpha\in A}$, 
and a subspace $W\subseteq C(X)$ such that there exists an isometric isomorphism 
$\phi\colon E\ra W$ with the following properties:

(a) For any $\alpha\in A$ and $f\in E_{\alpha,+,1}^{'}$ there exists $x_f\in X$ 
such that $f(e)=\phi(e)(x_f)$ for all $e\in E$.  

(b) If $f\colon E\to \RR$ with $f\in E^{'}_{\alpha,+}$ for some $\alpha\in A$, 
then there exists a BNN extension $\tl{f}\in C(X)^{'}_{\alpha,+}$ of $f$, 
where $E$ is identified with its image in $C(X)$. 
Moreover if 
the family of seminorms $\{ p_{\alpha} \}_{\alpha\in A}$ is saturated, then for 
any positive and continuous functional $f\colon E\to \RR$ there exists $\alpha\in A$ 
such that $f\in E^{'}_{\alpha,+}$ and consequently 
there exists a BNN extension $\tl{f}\in C(X)^{'}_{\alpha,+}$ of $f$. 
\end{theorem}

\begin{proof}
(a) We repeat Schaefer's construction in the proof of Theorem \ref{t:schaf} with the extra assumption 
on the existence of certain types of states. Consider a locally compact Hausdorff space $X$ and a map 
 $\phi\colon E\to C(X)$ as in the proof of Theorem \ref{t:schaf}. The map $\phi$ is correctly defined and linear, 
and is bipositive, 
as shown in the proof of Theorem \ref{t:schaf}. In order to see that $\phi$ is an injection, 
let us show that $\mathrm{ker}(\phi)=\{0\}$. For a contradiction, assume that there exists 
$x\neq 0$ in $E$ such that $0=\phi(x)(f,\alpha)=f(x)$ for all $(f,\alpha)\in X$. Since 
$E$ is Hausdorff separated, there exists $\alpha\in A$ such that $p_{\alpha}(x)\neq 0$. Then by the assumption there exists 
$f\in E_{\alpha,+,1}^{'}$ with $\vert f(x) \vert=p_{\alpha}(x)$, 
a contradiction. Thus, $\phi$ is an injection.  

In order to show that $\phi$ is an isometry, hence both $\phi$ and $\phi^{-1}$ are continuous, we use the same argument in the proof of 
Theorem \ref{t:schaf} that 
showed isometry on the positive elements of $E$, only this time due to the stronger extra assumption, the isometry 
works on all elements of $E$. Let $K_{\alpha}:=B_{\alpha}\times\{ \alpha \}$ 
for any $\alpha\in A$. Then each $K_{\alpha}$ is compact, the sets $\{ K_{\alpha} \}_{\alpha\in A}$ are pairwise disjoint 
with $X:=\cup_{\alpha\in A} K_{\alpha}$ 
and we consider the supremum on $K_{\alpha}$ seminorms $\{ \vert \cdot  \vert_{K_{\alpha}} \}_{\alpha\in A}$. 
Now given an arbitrary $\alpha\in A$ and $x\in E$ we have 
\begin{align*}
\vert \phi(x) \vert_{K_{\alpha}}&= \underset{(f,\alpha)\in K_{\alpha}}{\mathrm{sup}}\,\vert \phi(x)(f,\alpha) \vert \\
&=\underset{(f,\alpha)\in K_{\alpha}}{\mathrm{sup}}\,\vert f(x) \vert \\
&=\underset{0\leq f\leq p_{\alpha}}{\mathrm{sup}}\,\vert f(x) \vert \\
&=p_{\alpha}(x)
\end{align*} 
therefore all the properties of $\phi$ are established. Let $W:=\phi(E)$. 

Now given any $\alpha\in A$ and $f\in E_{\alpha,+,1}^{'}$, consider $x_f:=(f,\alpha)\in K_{\alpha}$. 
Then for any $e\in E$ we have $f(e)=\phi(e)(f,\alpha)=\phi(x_f)$ for all $e\in E$ and part (a) is shown. 

(b) Let us denote the induced functional on $W$, see Remark \ref {r:indfl}, again by $f$. We 
will now check the Bauer-Nachbin-Namioka Theorem, see Theorem \ref{t:baunami}, on $f$. 
Let us consider the normalized $f^{'}=f/\vert f \vert_{\alpha}$ so that $f^{'}\in  E^{'}_{\alpha,+,1}$. 
Given such $\alpha\in A$, consider the supremum on $K_{\alpha}$ seminorm, $\vert \cdot \vert_{\alpha}$. 
Assume that $e\in E$ and $l\in C(X)_+$ are such that $\vert \phi(e)+l \vert_{\alpha}\leq 1$, equivalently, 
$\underset{t\in K_{\alpha}}{\mathrm{sup}} \vert (\phi(e)+l)(t) \vert \leq 1$. Since $l(t)\geq 0$ for any $t\in X$, 
we have $\phi(e)(t)\leq 1$ for any $t\in K_{\alpha}$. 
Since $\phi(e)(t)$ for any 
$t\in K_{\alpha}$ is 
of the form $\phi(e)(g,\alpha)=g(e)$ for some $g\in  E^{'}_{\alpha,+,1}$,  
it follows that $f^{'}(e)\leq 1$. By the Bauer-Nachbin-Namioka Theorem, there exists a positive 
functional $\tl{f}\colon C(X)\ra \RR$ extending $f$ with $\vert \tl{f} \vert_{\alpha}=\vert f \vert_{\alpha}$. 
In particular, $\tl{f}\in C(X)^{'}_{\alpha,+}$. 

Finally, if $f\colon E\to \RR$ is any positive and continuous functional, assuming that $\{ p_{\alpha} \}_{\alpha\in A}$ is saturated, 
there exists $\alpha\in A$ such that $f\in E^{'}_{\alpha,+}$. Applying the result above 
there exists a BNN extension $\tl{f}\in C(X)^{'}_{\alpha,+}$ of $f$. 
\end{proof}

\begin{remark}\label{r:satur}
Notice that, compared to Theorem \ref{t:schaf}, Theorem \ref{t:posext} does not require the saturatedness assumption 
on the specified family of seminorms $\{ p_{\alpha}\}_{\alpha\in A}$ due to the stronger assumption on 
the existence of a certain type of state. 
\end{remark}

In the following we make a discussion of quotient spaces of a locally convex space 
by kernels of its seminorms, which will be useful for the rest of the article. 
For any locally convex space $(M, \{ q_{\beta} \}_{\beta\in B})$ 
with a specified family of continuous generating seminorms $\{ q_{\beta} \}_{\beta\in B}$, 
let us consider the quotient normed space $(M)_{q_{\beta}}$ for a fixed $\beta\in B$, 
with the natural norm given by $\Vert x_{q_{\beta}} \Vert_{q_{\beta}}:= q_{\beta}(x)$ for any $x\in M$. 
Let $f\in (M)_{q_{\beta}}^{'}$ be a functional with norm $\Vert f \Vert$. Define $g\colon M\ra \RR$ 
by $g(x):=f(x_{q_{\beta}})$. Notice that for any $x,y\in M$ with 
$q_{\beta}(x-y)=0$ we have $g(x)=g(y)$. It is easy to check that 
$g$ is linear, and we have 
$\vert g(x)\vert =\vert f(x_{q_{\beta}}) \vert \leq \Vert f \Vert \Vert x_{q_{\beta}} \Vert_{q_{\beta}}
 = \Vert f \Vert q_{\beta}(x)$ 
for all $x\in M$, therefore $g\in M_{\beta}^{'}$ and $\vert g \vert_{\beta}\leq \Vert f \Vert$. 
Similarly one obtains $\Vert f \Vert \leq \vert g \vert_{\beta}$, so $\vert g \vert_{\beta} = \Vert f \Vert$. 

Let now $(E, E_+, \{ p_{\alpha} \}_{\alpha\in A})$ be a calibrated normal ordered space with a closed cone $E_+$. 
For any fixed $\alpha\in A$, we consider the quotient normed ordered space $((E)_{p_{\alpha}}, \ol{(E_+)_{p_{\alpha}}}, \Vert \cdot \Vert_{p_{\alpha}})$ 
with the closure of the natural quotient cone $(E_+)_{p_{\alpha}}:=\{x_{\alpha} \mid x\in E_+  \}$. It is well 
known that, e.g. see Corollary 2.24 of \cite{Aliprantis} or Lemma 2.3 of \cite{Ay1}, 
the quotient space is then again a normal ordered space. Moreover, since any positive continuous 
functional $f\in ((E)_{p_{\alpha}})_{+}^{'}$ with respect to the cone $(E_+)_{p_{\alpha}}$ extends to the closure 
$\ol{(E_+)_{p_{\alpha}}}$, it is easy to see that the set $((E)_{p_{\alpha}})_{+}^{'}$ is the same with respect 
to both cones, so we can refer to it without ambiguity. Consider a positive functional 
$f\in ((E)_{p_{\alpha}})_{+}^{'}$ and define a functional $g\colon E\ra \RR$ by 
$g(x):=f(x_{p_{\alpha}})$. It is clear that $g$ is positive and by the paragraph before 
we have $g\in E^{'}_{\alpha,+}$ with $\vert g \vert_{\alpha}=\Vert f \Vert$. If, 
conversely, $f\in E^{'}_{\alpha,+}$ is a given positive 
functional, the quotient functional $f_{p_{\alpha}}(x_{p_{\alpha}}):=f(x)$ 
is well defined and positive on the space $((E)_{p_{\alpha}}, \ol{(E_+)_{p_{\alpha}}}, \Vert \cdot \Vert_{p_{\alpha}})$ 
and by the discussion above, for the corresponding operator norm we have 
$\Vert f_{p_{\alpha}} \Vert = \vert f \vert_{\alpha}$, in particular, 
$f_{p_{\alpha}} \in ((E)_{p_{\alpha}})_{+}^{'}$. We will 
use these observations in the rest of the article. 

In the following, given a calibrated ordered space $(C(X),C(X)_+,\{K_{\alpha}\}_{\alpha\in A})$ we will consider 
the quotient space $(C(X))_{p_{\alpha}}$ for a fixed seminorm $p_{\alpha}$, which, as before, 
is the seminorm of supremum on the compact set $K_{\alpha}$.  

We will 
show that $(C(X))_{p_{\alpha}}$ is isometrically isomorphic to the ordered Banach space 
$C(K_{\alpha})$ of continuous real valued functions on the compact set $K_{\alpha}$. 
Let us define the map $\xi\colon (C(X))_{p_{\alpha}}\ra C(K_{\alpha})$ by 
$\xi(x_{p_{\alpha}}):= x_{|{K_{\alpha}}}$, where $x_{|{K_{\alpha}}}$ 
is the restriction of the function $x$ to $K_{\alpha}$. Then $\xi$ 
is well defined, if $x_{p_{\alpha}}=y_{p_{\alpha}}$, then by definition 
$p_{\alpha}(x-y)=0$ and it follows that $x=y$ on the set $K_{\alpha}$. By reversing 
this argument one shows that $\xi$ is injective, and it is easy to see that 
$\xi$ is linear. Also, by the definition of $C(X)_+$ it is easy to see that $\xi$ is 
a positive map. In addition, we have 
\begin{align*}
\Vert x_{p_{\alpha}} \Vert_{p_{\alpha}}&=p_{\alpha}(x)
=\underset{t\in K_{\alpha}}{\mathrm{sup}} \vert x(t) \vert \\
&=\Vert x_{|{K_{\alpha}}} \Vert = \Vert \xi(x_{p_{\alpha}}) \Vert 
\end{align*}
for any $x\in C(X)$ and therefore $\xi$ is isometric. 

In order to show that the map $\xi$ is onto, we use a standart corollary of the Tietze 
extension theorem and the Urysohn lemma, see e.g section 4.5 of \cite{Folland}, 
and for convenience we provide a proof.  
Given $f\in C(K_{\alpha})$, we claim that there exists an extension 
$F\in C(X)$. By local compactness there exists an open set $U\supseteq K_{\alpha}$ 
such that $\ol{U}$ is compact. Since $\ol{U}$ is then compact and Hausdorff, it is a 
normal topological space, so the Tietze extension theorem and the Urysohn's lemma applies. 
By the Tietze extension theorem on $f$, there exists $\tl{f}\in C(\ol{U})$ extending $f$ 
and by the Urysohn's lemma applied with the closed sets $K_{\alpha}$ and $\ol{U}\setminus U$, 
there exists a continuous function $g\colon \ol{U}\ra [0,1]$ 
with $g_{|K_{\alpha}}=1$ and $g_{|\ol{U}\setminus U}=0$, where, WLOG 
we can assume $\ol{U}\setminus U \neq \emptyset$. Now define $F\colon X\ra \RR$ by 
\begin{equation*}
 F(X)\colon=\begin{cases} 
g(x)\tl{f}(x),\,\,&x\in \overline{U} \\ 0, &x \notin \ol{U}.
\end{cases} 
\end{equation*} 
Clearly $F$ is an extension of $f$ and it is not difficult to check that $F$ is continuous. 
Therefore, given any $f\in C(K_{\alpha})$ there exists $F\in C(X)$ such that 
$\xi(F_{p_{\alpha}})=f$, so $\xi$ is surjective. Moreover, it is clear that, 
if $f\geq 0$ on $K_{\alpha}$, then $F\geq 0$ on $X$, therefore, the inverse 
$\xi^{-1}\colon C(K_{\alpha}) \ra (C(X))_{p_{\alpha}}$ is a positive map as well, 
so $\xi$ is a bipositive map. Therefore, $(C(X))_{p_{\alpha}}$ is isometrically 
isomorphic to $ C(K_{\alpha})$. 

As a result, in the remainder of this section 
we will identify $(C(X))_{p_{\alpha}}$ with $C(K_{\alpha})$ and 
we will also identify functionals, in particular positive functionals 
on $(C(X))_{p_{\alpha}}$ with their induced 
functionals on $C(K_{\alpha})$, see Remark \ref{r:indfl}. 

\begin{remark}
In fact, $(C(X))_{p_{\alpha}}$  
is a real algebra as well by defining $x_{p_{\alpha}} y_{p_{\alpha}}:= (xy)_{p_{\alpha}}$. 
The fact that this multiplication is well defined can be seen using the fact 
that $C(X)$ is a real algebra and $p_{\alpha}$ multiplicative, i.e. $p_{\alpha}(xy)=p_{\alpha}(x)p_{\alpha}(y)$ 
for all $x,\in C(X)$. The map 
$\xi$ respects multiplication 
as well, since we have 
$\xi(x_{p_{\alpha}}y_{p_{\alpha}})=(xy)_{|{K_{\alpha}}}=x_{|{K_{\alpha}}}y_{|{K_{\alpha}}}
=\xi(x_{p_{\alpha}})\xi(y_{p_{\alpha}})$ for all $x_{p_{\alpha}},y_{p_{\alpha}}\in (\cC(X))_{p_{\alpha}}$, 
so $\xi$ is also a real algebra isomorphism, but we will not use that fact in this article. 
\end{remark}

The following theorem is a locally convex version of a combination of the Riesz-Markov-Kakutani Theorem and 
the Jordan decomposition for measures for the case of $C(X)$, where $X$ 
is a locally compact Hausdorff space. It will be needed in the proofs of the main theorems below. 

\begin{theorem}\label{t:rmlocal}
Let $X$ be a locally compact Hausdorff space and $(C(X),C(X)_+,\{ K_{\alpha} \}_{\alpha\in A})$ be a 
calibrated ordered space. Let $\alpha\in A$ be arbitrary and $u\in C(X)^{'}_{\alpha}$. Then 
there exists $v_1,v_2\in C(X)^{'}_{\alpha,+}$ such that $u=v_1-v_2$ and 
$\vert u \vert_{\alpha} = \vert v_1 \vert_{\alpha} + \vert v_2 \vert_{\alpha}$. 
\end{theorem}

\begin{proof}
By the discussion above we identify $(C(X))_{p_{\alpha}}$ with $C(K_{\alpha})$ and 
we will also identify functionals on $(C(X))_{p_{\alpha}}$ with their induced 
functionals on $C(K_{\alpha})$, see Remark \ref{r:indfl}. 

Let $u\in C(X)^{'}_{\alpha}$ be a functional. Consider the 
quotient functional $u_{p_{\alpha}}\colon (C(X))_{p_{\alpha}}\ra \RR$ 
defined by $u_{p_{\alpha}}(x_{p_{\alpha}}):=u(x)$ for any $x\in C(X)$. Note 
that we have $\vert u \vert_{\alpha}= \Vert u_{p_{\alpha}} \Vert$ by the above discussion  
and as a result we have $u_{p_{\alpha}} \in (C(X))^{'}_{p_{\alpha}}$. By the above discussion again 
we have $w:=u_{p_{\alpha}}\in C(K_{\alpha})^{'}$, since $(C(X))_{p_{\alpha}}=C(K_{\alpha})$. By a well known corollary of the 
Riesz-Markov-Kakutani Theorem, e.g. see the discussion on p.91 of \cite{Murphy}, 
there exist two positive linear functionals $w_1,w_2\colon C(K_{\alpha})\ra \RR$ 
with $w=w_1-w_2$ and $\Vert w \Vert = \Vert w_1 \Vert + \Vert w_2 \Vert $.  

Now let us define functionals $v_1,v_2\colon C(X) \ra \RR$ by 
$v_1(x):=w_1(x_{p_{\alpha}})$ and $v_2(x):=w_2(x_{p_{\alpha}})$. By the discussion 
before the theorem, we have $v_1,v_2\in C(X)^{'}_{\alpha,+}$ with $\vert v_1 \vert_{\alpha} = \Vert w_1 \Vert$ 
and $\vert v_2 \vert_{\alpha} = \Vert w_2 \Vert$. Note that for any $x, y\in C(X)$, 
if $x_{p_{\alpha}} = y_{p_{\alpha}}$, then $p_{\alpha}(x-y)=0$ and therefore 
$v_i(x-y)\leq p_{\alpha}(x-y)=0$, $i=1,2$ so $v(x)=v(y)$. Consequently we have 
\begin{align*}
u(x)&=u_{p_{\alpha}}(x_{p_{\alpha}}) \\
&=w_1(x_{p_{\alpha}}) - w_2(x_{p_{\alpha}}) \\
&=v_1(x)-v_2(x)
\end{align*} 
for all $x\in C(X)$. Finally we have 
\begin{align*}
\vert u \vert_{\alpha} &= \Vert u_{p_{\alpha}} \Vert \\
&=\Vert w_1 \Vert + \Vert w_2 \Vert \\
&=\vert v_1 \vert_{\alpha} + \vert v_2 \vert_{\alpha}
\end{align*} 
and the proof is finished. 
\end{proof}

Theorems \ref{t:isorepcx}, \ref{t:isorepsnfl}, \ref{t:isorepstr} and \ref{t:isorepquot} are the main theorems of this article. They all give 
equivalent characterizations of isometric representability of a calibrated ordered space $E$ on a 
space of continuous functions $C(X)$ where $X$ is a locally compact Hausdorff space, also see Remarks \ref{r:comchar} and 
\ref{r:refofisorep} below. 

\begin{theorem}\label{t:isorepcx}
Let $(E,E_+,\{ p_{\alpha}\}_{\alpha\in A})$ be a calibrated ordered space with a closed cone $E_+$. Then the following are equivalent. 

\begin{enumerate}
\item For any $\alpha\in A$ and  
$x_0\in E$ there exists a $p_{\alpha}$-state $f\in E_{\alpha,+,1}^{'}$ with $\vert f(x_0) \vert=p_{\alpha}(x_0)$. \label{isorep5}

\item There exists a locally compact Hausdorff space $X$ and a family of disjoint compact subsets $\{K_{\alpha}\}_{\alpha\in A}$ with 
$X=\bigcup_{\alpha\in A} K_{\alpha}$, such that any 
compact subset $K\subseteq X$ is contained in a finite union from the 
family $\{ K_{\alpha}\}_{\alpha\in A}$,  
and a subspace $W\subseteq C(X)$, where $W$ is isometrically isomorphic to $E$, 
such that if $f\in E^{'}_{\alpha,+}$ for some $\alpha\in A$, 
then there exists a BNN extension $\tl{f}\in C(X)^{'}_{\alpha,+}$. 
Moreover if 
the family of seminorms $\{ p_{\alpha} \}_{\alpha\in A}$ is saturated, then for 
any positive and continuous functional $f\colon E\to \RR$ there exists $\alpha\in A$ 
such that $f\in E^{'}_{\alpha,+}$ and consequently 
there exists a BNN extension $\tl{f}\in C(X)^{'}_{\alpha,+}$ of $f$. \label{isorep11}

\item There exists a locally compact Hausdorff space $X$ and a family of disjoint compact subsets $\{K_{\alpha}\}_{\alpha\in A}$ with 
$X=\bigcup_{\alpha\in A} K_{\alpha}$, such that any 
compact subset $K\subseteq X$ is contained in a finite union from the 
family $\{ K_{\alpha}\}_{\alpha\in A}$, 
and a subspace $W\subseteq C(X)$, where $W$ is isometrically isomorphic to $E$, 
such that for any functional $f\in E^{'}_{\alpha,+,1} $ for some 
	$\alpha\in A$, there exists $x_{f}\in K_{\alpha}$ such that $\phi(u)(x_{f})=f(u)$ for any $u\in E$. \label{isorep0}

\item There exists a locally compact Hausdorff space $X$ and a family of disjoint compact subsets $\{K_{\alpha}\}_{\alpha\in A}$ 
closed under finite unions with 
$X=\bigcup_{\alpha\in A} K_{\alpha}$, such that any 
compact subset $K\subseteq X$ is contained in a finite union from the 
family $\{ K_{\alpha}\}_{\alpha\in A}$,  
  and a subspace $W\subseteq C(X)$ such that $W$ is 
isometrically isomorphic to $E$. \label{isorep1}

\item There exists a Hausdorff topological space $X$ and a family of compact subsets $\{ K_{\alpha} \}_{\alpha\in A}$ closed under finite unions and 
$\bigcup_{\alpha\in A} K_{\alpha} = X$ such that there exists an isometric isomorphism $\phi\colon E\ra W$ where $W\subseteq C(X)$ is a subspace. \label{isorep12}

\end{enumerate}
\end{theorem}

\begin{proof}
``$(\ref{isorep5})\Ra(\ref{isorep11})$": See part (b) of Theorem \ref{t:posext}. 

``$(\ref{isorep5})\Ra(\ref{isorep0})$": See part (a) of Theorem \ref{t:posext}. 

``$(\ref{isorep11})\Ra(\ref{isorep1})$": Clear. 

``$(\ref{isorep0})\Ra(\ref{isorep1})$": Clear. 

``$(\ref{isorep1})\Ra(\ref{isorep12})$": Clear. 

``$(\ref{isorep12})\Ra(\ref{isorep5})$": Fix any $x_0\in E$ and $\alpha\in A$. Since $K_{\alpha}$ is compact 
and $\phi(x_0)\in W$ is continuous on $K_{\alpha}$, there exists $y_0\in K_{\alpha}$ such that 
$\phi(x_0)(y_0)$ or $-\phi(x_0)(y_0)$ realizes $\vert \phi(x_0) \vert_{\alpha}$. Consider the point evaluation functional 
$f_{y_0}\colon W \ra \RR$ given by $f_{y_0}(\phi(e)):=\phi(e)(y_0)$ for any $e\in E$. It is easy to check that 
$f_{y_0}$ is a linear positive functional. Notice that we have $\vert f_{y_0}(\phi(e)) \vert \leq \vert \phi(e) \vert_{\alpha}$ 
for any $e\in E$ and we also have 
$\vert \phi(x_0)(y_0) \vert \leq \vert \phi(x_0) \vert_{\alpha}=\vert f_{y_0}(\phi(x_0)) \vert$. Therefore 
$f_{y_0}\in W^{'}_{\alpha,+,1}$ and by Remark \ref{r:indfl} 
the existence of a functional with the required properties is shown.  
\end{proof}

The next theorem gives characterizations in terms of the generating family of seminorms and positive functionals of the space $E$.  

\begin{theorem}\label{t:isorepsnfl}

Let $(E,E_+,\{ p_{\alpha}\}_{\alpha\in A})$ be a calibrated ordered space with a 
closed cone $E_+$. Then the following are equivalent. 

\begin{enumerate}
\item For any $\alpha\in A$ and  
$x_0\in E$ there exists a $p_{\alpha}$-state $f\in E_{\alpha,+,1}^{'}$ with $\vert f(x_0) \vert=p_{\alpha}(x_0)$. \label{isorep51}

\item For any $\alpha\in A$ and 
any $x\in E$ at least one of the following properties holds: $p_{\alpha}(x+l)\geq p_{\alpha}(x)$ for all $l\in E_+$ and/or $p_{\alpha}(x-l)\geq p_{\alpha}(x)$ for all $l\in E_+$, i.e. all elements of $E$ are either semi positive and/or semi negative. \label{isorep3}

\item The corresponding unit balls 
$E_{\alpha,1}$ to each $\{p_{\alpha}\}_{\alpha\in A}$ are full; i.e. for any $\alpha\in A$ and $a\leq x \leq b$ with $p_{\alpha}(a),\,p_{\alpha}(b)\leq 1$ we have $p_{\alpha}(x)\leq 1$. \label{isorep2}

\item For all $x,y,z\in E$ with $x\leq y \leq z$ and any $\alpha\in A$ we have $p_{\alpha}(y)\leq \mathrm{max}\{ p_{\alpha}(x), p_{\alpha}(z)\}$. \label{isorep9} 

\item For any $\alpha\in A$ and $u\in E_{\alpha}^{'}$, 
there exists $v_1,v_2 \in E_{\alpha,+}^{'}$ such that $u=v_1-v_2$ and $\vert u \vert_{\alpha}=\vert v_1 \vert_{\alpha}+\vert v_2 \vert_{\alpha}$. \label{isorep4}

\item For any $\alpha\in A$ and $x_0\in E$ the equality $\underset{f\in E^{'}_{\alpha,+,1}}{\mathrm{sup}}\vert f(x_0) \vert = p_{\alpha}(x_0)$ holds. \label{isorep19}

\end{enumerate}
\end{theorem}

\begin{proof}
``$(\ref{isorep51})\Ra(\ref{isorep3})$": By Theorem \ref{t:isorepcx}, $(\ref{isorep51})$ is equivalent with $(\ref{isorep12})$ of 
Theorem \ref{t:isorepcx}, so we show that $(\ref{isorep12})$ of Theorem \ref{t:isorepcx} implies $(\ref{isorep3})$. 
Fix any $\alpha\in A$ and $x\in E$. Assume for a while that 
$\mathrm{sup}\{\phi(x)(t)\mid t\in K_{\alpha}\,\mathrm{and}\, \phi(x)(t)\geq 0\} 
\geq \mathrm{sup}\{-\phi(x)(t)\mid t\in K_{\alpha}\,\mathrm{and}\, \phi(x)(t)\leq 0\}$, where 
WLOG both sets are nonempty. 
Therefore $\vert \phi(x) \vert_{\alpha} = \mathrm{sup}\{\phi(x)(t)\mid t\in K_{\alpha}\,\mathrm{and}\, \phi(x)(t)\geq 0\}$.  
If $l\in E_+$ is arbitrary, it is then easy to see that 
$\vert \phi(x)+\phi(l) \vert_{\alpha} = \mathrm{sup}\{ (\phi(x+ l)(t)\mid t\in K_{\alpha}\,\mathrm{and}\, (\phi(x+l)(t)\geq 0 \}$ 
and therefore $p_{\alpha}(x + l) = \vert \phi(x + l) \vert_{\alpha} \geq \vert \phi(x) \vert_{\alpha} = p_{\alpha}(x)$, 
so $x$ is semi-positive. 

If, on the other hand, 
$\mathrm{sup}\{-\phi(x)(t)\mid t\in K_{\alpha}\,\mathrm{and}\, \phi(x)(t)\leq 0\} 
\geq \mathrm{sup}\{\phi(x)(t)\mid t\in K_{\alpha}\,\mathrm{and}\, \phi(x)(t)\geq 0\}$, 
then a similar consideration shows that $p_{\alpha}(x - l) \geq p_{\alpha}(x)$ for all $l\in E_+$, 
so $x$ is semi-negative, and the implication is proven. 

``$(\ref{isorep3})\Ra(\ref{isorep2})$": Take any $\alpha\in A$ and assume that $a\leq x \leq b$ with 
$p_{\alpha}(a), p_{\alpha}(b) \leq 1$. If $x$ has $p_{\alpha}(x+l)\geq p_{\alpha}(x)$ for all $l\in E_+$ then 
$p_{\alpha}(b)\geq p_{\alpha}(x)$, so $p_{\alpha}(x) \leq 1$. The argument is similar if $x$ has $p_{\alpha}(x-l)\geq p_{\alpha}(x)$ 
for all $l\in E_+$, so $p_{\alpha}(x) \leq 1$ and the implication is proven. 

``$(\ref{isorep2})\Ra(\ref{isorep51})$": Let $x_0\in E$ and $\alpha\in A$ be arbitrary. Define a one dimensional subspace 
$X:=\{\lambda x_0 \mid \lambda\in \RR \}$. Define a linear functional $f\colon X\ra \RR$ by $f(\lambda x_0):=\lambda p_{\alpha}(x_0)$. 
If $\lambda \leq 0$ and $\lambda x_0 \in X \cap (A_{\alpha}-E_+)$, where 
$A_{\alpha}:=\{x\in E \mid p_{\alpha}(x)\leq 1 \}$, then $f(\lambda x_0)= \lambda p_{\alpha}(x_0)\leq 0 < 1$. Now 
let us denote ${}_{+}X:=\{\lambda x_0 \mid \lambda > 0 \}$. Assume that 
$x\in {}_{+}X \cap (A_{\alpha}-E_+)$, so $x=\lambda_0 x_0$ for some $\lambda_0 >0$, 
then $p_{\alpha}(x+l)\leq 1$ for some $l\in E_+$. Notice that in this case for any $0\leq \eta \leq \lambda_0$ 
we have $p_{\alpha}(\eta x_0 + l)\leq 1$ for some $l\in E_+$. 

We now consider two cases: \smallskip 

\emph{Case I:} If for all $x\in {}_{+}X$ such that $p_{\alpha}(x+l)\leq 1$ for some $l\in E_+$ there exists 
$m\in E_+$ such that $p_{\alpha}(x-m)\leq 1$, since we have $x-m \leq x \leq x+l$ and $p_{\alpha}(x+l), p_{\alpha}(x-m) \leq 1$, 
it follows that we have $p_{\alpha}(x)\leq 1$, equivalently, $f(x)\leq 1$. By the Bauer-Nachbin-Namioka Theorem, see Theorem \ref{t:baunami}, 
there exists a positive functional $f\colon E\ra \RR$ with $f(x_0)=p_{\alpha}(x_0)$ and $\vert f(x) \vert \leq p_{\alpha}(x)$ for all $x\in E$. 
\smallskip  

\emph{Case II:} If the condition in \emph{Case I} does not hold, then there exists $x\in {}_{+}X$ such that $p_{\alpha}(x+l)\leq 1$ 
for some $l\in E_+$ with the property that $p_{\alpha}(x-m) > 1$ for all $m\in E_+$. As a result we obtain $x-m \leq x \leq x+l$ 
with $p_{\alpha}(x-m) >1$ and $p_{\alpha}(x+l) \leq 1$ for some $l\in E_+$ and all $m\in E_+$. Dividing all sides by 
$p_{\alpha}(x-m)$ in the vector inequality we obtain 
\begin{equation*}
\frac{x-m}{p_{\alpha}(x-m)} \leq \frac{x}{p_{\alpha}(x-m)} \leq \frac{x+l}{p_{\alpha}(x-m)}.
\end{equation*}

Taking $p_{\alpha}$ of the first and the third terms, we obtain
\begin{equation*}
p_{\alpha}\left(\frac{x-m}{p_{\alpha}(x-m)}\right),p_{\alpha}\left(\frac{x+l}{p_{\alpha}(x-m)}\right)\leq 1
\end{equation*}
so we must have $\displaystyle{p_{\alpha}\left(\frac{x}{p_{\alpha}(x-m)}\right)\leq 1}$, 
equivalently, $p_{\alpha}(x)\leq p_{\alpha}(x-m)$ for all $m\in E_+$. From this an easy argument 
shows that $p_{\alpha}(x_0)\leq p_{\alpha}(x_0-m)$ for all $m\in E_+$. Thus by Corollary \ref{c:bau2} 
there exists a positive functional $f\colon E\ra \RR$ with $f(x_0)=-p_{\alpha}(x_0)$ and 
$\vert f(x) \vert \leq p_{\alpha}(x)$ for all $x\in E$. 

Combining \emph{Case I} and \emph{Case II} a proof of the implication is completed. 

``$(\ref{isorep51})\Ra(\ref{isorep4})$": We show that $(\ref{isorep1})$ of Theorem \ref{t:isorepcx} 
implies $(\ref{isorep4})$ as it is equivalent with $(\ref{isorep51})$ by that theorem. Given $u\in E^{'}_{\alpha}$ for some $\alpha\in A$, 
by the Hahn-Banach Theorem there exists an extension $\tl{u}\colon C(X)\ra \RR$ 
with $\tl{u}\in C(X)^{'}_{\alpha}$ and $\vert \tl{u} \vert_{\alpha} = \vert u \vert_{\alpha}$. 
Then by Theorem \ref{t:rmlocal} there exists $v_1,v_2\in C(X)^{'}_{\alpha,+}$ such that 
$\tl{u}=v_1-v_2$ and $\vert \tl{u} \vert_{\alpha}=\vert v_1 \vert_{\alpha}+\vert v_2 \vert_{\alpha}$. 
Now restricting $v_1,v_2$ to $E$, we obtain linear and positive functionals $v_1^{'},v_2^{'}: E\ra \RR$ 
with $u=v_1^{'}-v_2^{'}$. 
Clearly $\vert v_1^{'} \vert_{\alpha} \leq \vert v_1 \vert_{\alpha}$ and 
$\vert v_2^{'} \vert_{\alpha} \leq \vert v_2 \vert_{\alpha}$. Finally we have 
\begin{align*}
\vert v_1^{'} \vert_{\alpha}+\vert v_2^{'} \vert_{\alpha} &\geq \vert u \vert_{\alpha} = \vert \tl{u} \vert_{\alpha} \\
&=\vert v_1 \vert_{\alpha} + \vert v_2 \vert_{\alpha} \\
&\geq \vert v_1^{'} \vert_{\alpha} + \vert v_2^{'} \vert_{\alpha}
\end{align*}
therefore $\vert u \vert_{\alpha} = \vert v_1^{'} \vert_{\alpha} + \vert v_2^{'} \vert_{\alpha}$ holds 
and the implication is proven.

``$(\ref{isorep4})\Ra(\ref{isorep51})$": We follow the same steps as in the proof of Theorem 4.2 in \cite{Riedl}, adopted 
to the current case. WLOG we can assume that $p_{\alpha}(x_0)\neq 0$ and then 
there exists a functional $u\in E^{'}_{\alpha,1}$ such that $u(x_0)=p_{\alpha}(x_0)$. 
By the assumption, there exists two positive linear functionals $v_1,v_2\in E^{'}_{\alpha,+}$ such that $u=v_1-v_2$ 
and $\vert v_1 \vert_{\alpha} + \vert v_2 \vert_{\alpha} = \vert u \vert_{\alpha} = 1$, therefore $v_1,v_2\in E^{'}_{\alpha,+,1}$. 
Now we have 
\begin{align*}
p_{\alpha}(x_0)&=u(x_0)=v_1(x_0) - v_2(x_0) \\
&= v_1(x_0) + v_2(-x_0) \\
&\leq (\vert v_1 \vert_{\alpha} + \vert v_2 \vert_{\alpha}) p_{\alpha}(x_0) \\
&=p_{\alpha}(x_0)
\end{align*}
from which it follows that $v_1(x_0)=\vert v_1 \vert_{\alpha} p_{\alpha}(x_0)$ and 
$-v_2(x_0)=\vert v_2 \vert_{\alpha} p_{\alpha}(x_0)$. Since both $v_1(x_0)$ and $v_2(x_0)$ 
cannot be zero, at least one of $\vert v_1 \vert_{\alpha}$ and $\vert v_2 \vert_{\alpha}$ 
is different then zero, assume WLOG $\vert v_1 \vert_{\alpha}\neq 0$. Then 
$f=v_1/\vert v_1 \vert_{\alpha} \in E^{'}_{\alpha,+,1}$ and $f(x_0)=p_{\alpha}(x_0)$, 
so $f$ is a functional with required properties. 

``$(\ref{isorep51})\Ra(\ref{isorep19})$": This is clear. Notice that, it is a similar corollary 
with that of Corollary \ref{c:bauay2}, only this time the equality holds for all $x\in E$ 
due to the stronger assumption. 

``$(\ref{isorep19})\Ra(\ref{isorep51})$": We show that $(\ref{isorep19})$ implies $(\ref{isorep0})$ of Theorem \ref{t:isorepcx}, 
equivalent to $(\ref{isorep51})$ by that theorem. We repeat the proof of part (a) of Theorem \ref{t:posext}. 
Only showing that the map $\phi$ is injective is slightly different: 
Let the space $X$ and the map $\phi\colon E\to C(X)$ be as in the proof of Theorem \ref{t:posext}. 
In order to see that $\phi$ is an injection, 
assume by contradiction that there exists 
$x\neq 0$ in $E$ such that $0=\phi(x)(f,\alpha)=f(x)$ for all $(f,\alpha)\in X$. Since 
$E$ is Hausdorff, there exists $\alpha\in A$ such that $p_{\alpha}(x)\neq 0$. 
By the assumption, $\underset{f\in E^{'}_{\alpha,+,1}}{\mathrm{sup}}\vert f(x) \vert = p_{\alpha}(x)$, 
so there exists 
$f\in E_{\alpha,+,1}^{'}$ with $\vert f(x) \vert > 0$, 
a contradiction. Thus, $\phi$ is an injection. 

``$(\ref{isorep2})\Ra(\ref{isorep9})$": Let $\alpha\in A$ be arbitrary and $x,y,z\in E$ be such that 
$x\leq y \leq z$. Let $A:=\mathrm{max}\{ p_{\alpha}(x), p_{\alpha}(z)\}$. Assume for a while that 
$A>0$. Then by dividing all sides by $A$ we have $\displaystyle{\frac{x}{A}\leq \frac{y}{A} \leq \frac{z}{A}}$ 
with $p_{\alpha}(x/A),p_{\alpha}(z/A)\leq 1$. Therefore $p_{\alpha}(y/A)\leq 1$, equivalently, 
$p_{\alpha}(y)\leq A$. 

In the case $A=0$, replacing $A$ with $A+\epsilon$ for any $\epsilon >0$ and a standart argument 
shows that $p_{\alpha}(y)=0$, so the implication is shown in this case as well. 

``$(\ref{isorep9})\Ra(\ref{isorep2})$": Let $\alpha\in A$ be arbitrary and $x,y,z\in E$ be such that 
$x\leq y \leq z$ with $p_{\alpha}(x), p_{\alpha}(z)\leq 1$. Since $p_{\alpha}(y)\leq \mathrm{max}\{ p_{\alpha}(x), p_{\alpha}(z) \}\leq 1$, 
the implication is proven. 
\end{proof}

\begin{remark}
Notice that, in the preceeding theorem the equivalence of $(\ref{isorep51})$ and $(\ref{isorep3})$ is also 
contained in Corollary \ref{c:bau3}. 
\end{remark}

The following theorem gives characterizations in terms of isometrically isomorphically representability as a subspace of 
different structures; namely, on locally $C^*$-algebras and locally Riesz spaces.  

\begin{theorem}\label{t:isorepstr}
Let $(E,E_+,\{ p_{\alpha}\}_{\alpha\in A})$ be a calibrated ordered space with a 
closed cone $E_+$. Then the following are equivalent. 

\begin{enumerate}
\item For any $\alpha\in A$ and  
$x_0\in E$ there exists a $p_{\alpha}$-state $f\in E_{\alpha,+,1}^{'}$ with $\vert f(x_0) \vert=p_{\alpha}(x_0)$. \label{isorep52}

\item There exists a commutative locally $C^*$-algebra $(\cA, \{p_{\alpha}\}_{\alpha\in A})$ and a real vector subspace $W\subseteq\cA^{h}$ such that $E$ is 
isometrically isomorphic to $W$. \label{isorep18}

\item There exists a locally $C^*$-algebra $(\cA, \{p_{\alpha}\}_{\alpha\in A})$ and a real vector subspace $W\subseteq\cA^{h}$ such that $E$ is 
isometrically isomorphic to $W$. \label{isorep7}

\item There exists a locally Riesz space $R$ and a family $\{ q_{\alpha} \}_{\alpha\in A}$ of continuous generating M-Riesz seminorms 
on $R$ and a subspace $P\subseteq R$ such that $P$ is isometrically isomorphic to $E$. \label{isorep10}

\item There exists an L-Archimedean local order unit space 
$(G,\{ C_{\alpha} \}_{\alpha\in A})$ with a local order unit $e$ and a subspace $H\subseteq G$ such that 
$E$ is isometrically isomorphic to $H$. \label{isorep60}
\end{enumerate}
\end{theorem}

\begin{proof}
``$(\ref{isorep52})\Ra(\ref{isorep18})$": We show that $(\ref{isorep1})$ of Theorem \ref{t:isorepcx} implies $(\ref{isorep18})$ 
as it is equivalent to $(\ref{isorep52})$ by that theorem. Let $\cA:=\cC(X)=C(X)+iC(X)$ 
with topology given by the seminorms supremum on the compact sets $\{ K_{\alpha}\}_{\alpha\in A}$. 
As discussed before, $\cA$ is a locally $C^*$-algebra with $(\cC(X))^{h}=C(X)$ and therefore 
the implication is proven. 

``$(\ref{isorep18})\Ra(\ref{isorep7})$": This is clear. 

``$(\ref{isorep7})\Ra(\ref{isorep52})$": Fix any $C^*$-seminorm $p_{\alpha}$ from the given $C^*$-calibration $\{ p_{\alpha} \}_{\alpha\in A}$ and 
any $x_0\in W$. Since $W\subseteq \cA^{h}$, note that $x_0$ is selfadjoint, hence normal. Therefore the quotient element 
${x_{0}}_{p_{\alpha}}\in \cA_{p_{\alpha}}$ also is normal. By the fact that the quotient 
$\cA_{p_{\alpha}}$ is a $C^*$-algebra, see \cite{Apostol} and 
a well known fact from the theory of $C^*$-algebras, e.g. see \cite{Murphy}, p.90, 
there exists a $C^*$-algebra state $g\colon \cA_{p_{\alpha}}\ra \CC$ such that 
$g({x_{0}}_{p_{\alpha}})=\Vert {x_{0}}_{p_{\alpha}} \Vert_{p_{\alpha}}=p_{\alpha}(x_0)$. 
Now define a functional $f\colon \cA\ra \CC$ by $f(x):= g(x_{p_{\alpha}})$. It is a routine check to see 
that $f\in \cA^{'}_{\alpha, +, 1}$ with $f(x_0)=p_{\alpha}(x_0)$. Restricting $f$ to $W$, 
we obtain a functional $f_1\colon W\ra \CC$. Since $f$ is in particular a positive functional, 
in fact $f_1\colon W\ra \RR$, and again it is routine to check that $f_1\in W^{'}_{\alpha, +, 1}$ 
with $f_1(x_0)=p_{\alpha}(x_0)$, 
so the existence of a required functional is shown. 

``$(\ref{isorep52})\Ra(\ref{isorep10})$": We show that $(\ref{isorep12})$ of Theorem \ref{t:isorepcx} implies 
$(\ref{isorep10})$ as it is equivalent to $(\ref{isorep52})$ by that theorem. By Example \ref{r:cxriesz} $C(X)$ is a locally Riesz space 
with topology generated by M-seminorms, so the implication holds. 

``$(\ref{isorep10})\Ra(\ref{isorep52})$": Let $\alpha\in A$ be arbitrary and note that since $q_{\alpha}$ is in particular 
an M-Riesz seminorm, $q_{\alpha}$ is increasing. Let $x_0\in R$ be arbitrary and assume for a while that, 
$q_{\alpha}(x_{0}^+)\vee q_{\alpha}(x_{0}^-)=q_{\alpha}(x_{0}^+)$, 
so that we have 
$q_{\alpha} (x_0 )=q_{\alpha}(\vert x_0 \vert)=q_{\alpha}(x_{0}^+\vee x_{0}^-)=q_{\alpha}(x_{0}^+)\vee q_{\alpha}(x_{0}^-)=q_{\alpha}(x_{0}^+)$. 
Now given any $x_0\in R$, by Corollary \ref{c:bauay}, there exists a positive functional 
$f\colon R\ra \RR$ such that $f(x_{0}^{+})=q_{\alpha}(x_{0}^+)$ and $\vert f(x) \vert\leq q_{\alpha}(x)$ for all $x\in R$. We now have 
\begin{align*}
q_{\alpha}(\vert x_0 \vert) \geq f(\vert x_0 \vert) &=f(x_0^+ + x_0^- ) \\
&= f(x_0^+ ) + f(x_0^- ) \\
&\geq q_{\alpha}(x_0^+)
\end{align*}
from which it follows that $f(x_0^-)=0$. We also have
\begin{align*}
q_{\alpha}( x_0 ) \geq f( x_0 ) &=f(x_0^+ - x_0^- ) \\
&= f(x_0^+ ) - f(x_0^- ) \\
&= f(x_0^+ ) 
\end{align*}
and therefore $f(x_0)=q_{\alpha}(x_0)$. Consequently, in this case, $(\ref{isorep5})$ holds. 

If, on the other hand, $q_{\alpha}(x_{0}^+)\vee q_{\alpha}(x_{0}^-)=q_{\alpha}(x_{0}^-)$ holds, so that we have 
$q_{\alpha}(x_0)=q_{\alpha}(\vert x_0 \vert)=q_{\alpha}(x_{0}^-)$, by arguing similarly to above, we first 
pick a positive functional 
$f\colon R\ra \RR$ such that $f(x_{0}^{-})=q_{\alpha}(x_{0}^-)$ and $\vert f(x) \vert\leq q_{\alpha}(x)$ for all $x\in R$, by Corollary 
\ref{c:bauay}. Similar calculation gives $f(x_0^+)=0$ and therefore $q_{\alpha}(x_0)\geq -f(x_0)=-f(x_0^+)+f(x_0^-)$. It follows that 
$q_{\alpha}(x_0)=f(x_0)$ and $(\ref{isorep5})$ is shown to hold in this case as well.  

``$(\ref{isorep52})\Ra(\ref{isorep60})$": We show that $(\ref{isorep12})$ of Theorem \ref{t:isorepcx} implies 
$(\ref{isorep60})$ as it is equivalent to $(\ref{isorep52})$ by that theorem. Since by Example \ref{r:cxlocou} $C(X)$ 
is an L-Archimedean local order unit space with an order unit the constant function $1$, the implication holds. 

``$(\ref{isorep60})\Ra(\ref{isorep52})$": We show that $(\ref{isorep60})$ implies $(\ref{isorep2})$ of Theorem \ref{t:isorepsnfl} 
as it is equivalent to $(\ref{isorep52})$ by that theorem; namely, we show that the unit balls $E_{\alpha,1}$, 
where we consider the order unit space as a calibrated ordered space by Example \ref{r:cxlocou}, are full for each $\alpha\in A$. 
Assume that $x\leq y\leq z$ with $p_{\alpha}(x),\, p_{\alpha}(z) \leq 1$. Then we have 
$-e\leq_{\alpha} x\leq y \leq z \leq_{\alpha} e$ and therefore $p_{\alpha}(y)\leq 1$, 
so $E_{\alpha,1}$ is full for each $\alpha\in A$ and the implication is shown. 
\end{proof}

We now state and prove equivalent characterizations in terms quotient spaces, see the discussion after 
Theorem \ref{t:posext}.  

\begin{theorem}\label{t:isorepquot}
Let $(E,E_+,\{ p_{\alpha}\}_{\alpha\in A})$ be a calibrated ordered space with a 
closed cone $E_+$. Assume also that the calibrated ordered space $(E,E_+,\{ p_{\alpha}\}_{\alpha\in A})$ is normal. 
Then the following are equivalent. 
\begin{enumerate}
\item For any $\alpha\in A$ and  
$x_0\in E$ there exists a $p_{\alpha}$-state $f\in E_{\alpha,+,1}^{'}$ with $\vert f(x_0) \vert=p_{\alpha}(x_0)$. \label{isorep53}

\item For any $\alpha\in A$ there exists a compact Hausdorff space $X_{\alpha}$ and 
a subspace of $W_{\alpha}$ of the space $C(X_{\alpha})$ such that the quotient space 
$(E)_{p_{\alpha}}$ is isometrically isomorphic to $W_{\alpha}$.  \label{isorep15}

\item For any $\alpha\in A$ there exists a Riesz space $R$ with topology given by an M-Riesz norm 
and a subspace $P\subseteq R$ such that $P$ is isometrically isomorphic to the quotient normed ordered space $(E)_{p_{\alpha}}$. \label{isorep16}

\item For any $\alpha\in A$ there exists an Archimedean order unit space $(F,F_+,f)$ and a subspace 
$K\subseteq F$ such that $K$ is isometrically isomorphic to the quotient normed ordered space $(E)_{p_{\alpha}}$. \label{isorep20}

\item For any $\alpha\in A$ there exists an almost Archimedean order unit space $(F,F_+,f)$ and a subspace 
$K\subseteq F$ such that $K$ is isometrically isomorphic to the quotient normed ordered space $(E)_{p_{\alpha}}$. \label{isorep17}

\end{enumerate}

\end{theorem}

\begin{proof}
``$(\ref{isorep53})\Ra(\ref{isorep15})$": We prove that $(\ref{isorep1})$ of Theorem \ref{t:isorepcx} 
implies $(\ref{isorep15})$ as it is equivalent to $(\ref{isorep53})$ by that theorem. Since $E$ is isometrically isomorphic to $W$, 
$(E)_{p_{\alpha}}$ is isometrically isomorphic to $(W)_{p_{\alpha}}\subseteq (C(X))_{p_{\alpha}}$.  
By the discussion preceeding Theorem \ref{t:rmlocal} we have $(C(X))_{p_{\alpha}}$ 
is isometrically isomorphic to $C(K_{\alpha})$, hence letting $X_{\alpha}:=K_{\alpha}$, 
it follows that $(E)_{p_{\alpha}}$ is isometrically isomorphic to a subspace of $C(X_{\alpha})$. 

``$(\ref{isorep15})\Ra(\ref{isorep53})$": Using $(\ref{isorep1})\Ra(\ref{isorep5})$ of Theorem \ref{t:isorepcx} for the 
special case of normed spaces, given any ${x_0}_{p_{\alpha}}$ there exists a positive continuous functional 
$f_{p_{\alpha}}\in ((E)_{p_{\alpha}})^{'}_{+,1}$ such that $\vert f_{p_{\alpha}}({x_0}_{p_{\alpha}})\vert=\Vert {x_0}_{p_{\alpha}} \Vert_{p_{\alpha}}$. 
Now define a functional $f\colon E\ra \RR$ by $f(x):=f_{p_{\alpha}}(x_{p_{\alpha}})$. By similar arguments as in the 
proof of Theorem \ref{t:rmlocal} it follows that $f$ is a linear and positive functional, in fact, 
$f\in E_{\alpha,+,1}^{'}$ with $\vert f(x_0) \vert = p_{\alpha}(x_0)$ and we are done.

``$(\ref{isorep15})\Ra(\ref{isorep16})$": Since $C(X_{\alpha})$ is a Riesz space with topology given 
by an M-Riesz norm by Example \ref{r:cxriesz}, the implication holds. 

``$(\ref{isorep16})\Ra(\ref{isorep15})$": By using the equivalence of $(\ref{isorep10})$ of Theorem \ref{t:isorepstr} and 
$(\ref{isorep1})$ of Theorem \ref{t:isorepcx} for the special case of normed ordered 
spaces, the implication follows. 

``$(\ref{isorep15})\Ra(\ref{isorep20})$": This is clear since $C(X_{\alpha})$ is an Archimedean order unit space. 

``$(\ref{isorep20})\Ra(\ref{isorep17})$": Clear, as any Archimedean order unit space is an almost Archimedean order unit space. 

``$(\ref{isorep17})\Ra(\ref{isorep15})$": Let us denote the norm generated by the order unit by $\Vert \cdot \Vert$, 
so $\Vert x \Vert=\underset{\lambda\geq 0}{\mathrm{inf}}\{ -\lambda f \leq x \leq \lambda f \}$. In order 
to show that the unit ball of $\Vert \cdot \Vert$ is full, assume that $x\leq y\leq z$ with 
$\Vert x \Vert,\, \Vert z \Vert \leq 1$. Then we have $-e\leq x \leq e$ and $-e\leq z \leq e$, 
so we obtain $-e\leq x\leq y \leq z \leq e$ and consequently $\Vert y \Vert \leq 1$, 
so the unit ball of $\Vert \cdot \Vert$ is full. Hence by the equivalence of $(\ref{isorep2})$ of 
Theorem \ref{t:isorepsnfl} and $(\ref{isorep1})$ of Theorem \ref{t:isorepcx}  
for the special case of normed ordered spaces, $(F,F_+,\Vert \cdot \Vert)$ is isometrically 
isomorphic to a subspace $W$ of $C(X_{\alpha})$ for a compact Hausdorff space $X_{\alpha}$. 
It follows that $E_{p_{\alpha}}$ is isometrically isomorphic to a subspace of 
$C(X_{\alpha})$ as well and the implication is proven. 

\end{proof}

\begin{remark}\label{r:comchar}
All of the theorems \ref{t:isorepcx}, \ref{t:isorepsnfl}, \ref{t:isorepstr} and \ref{t:isorepquot} share the common characterization that, 
for any $\alpha\in A$ and  
$x_0\in E$ there exists a $p_{\alpha}$-state $f\in E_{\alpha,+,1}^{'}$ with $\vert f(x_0) \vert=p_{\alpha}(x_0)$ on which 
most of the proofs are based. In particular, all of the characterizations appearing in theorems 
\ref{t:isorepcx}, \ref{t:isorepsnfl}, \ref{t:isorepstr} are logically equivalent and if in addition 
the space $E$ is normal, they are logically equivalent with the characterizations in Theorem \ref{t:isorepquot} 
as well. 
\end{remark}

\begin{remark}\label{r:onto}
It is interesting to consider the natural question when the calibrated ordered space $E$ is isomorphically 
isomorphic to a space $C(X)$ where $X$ is locally compact and Hausdorff. In the normed case, some characterizations 
are obtained in Theorem 5 of \cite{Riedl} and for AM (Abstract M)-spaces, in \cite{Kakutani}. 
In this article we do not go in this direction. 
\end{remark}

By considering the special case of normed ordered spaces, we obtain the following corollary:

\begin{theorem}\label{t:isorep2}
Let $(E,E_+,\Vert \cdot \Vert)$ be a normed ordered space with a 
closed cone $E_+$. Then the following are equivalent. 

\begin{enumerate}
\item There exists a compact Hausdorff space $X$ 
and a subspace $W\subseteq C(X)$, where $W$ is isometrically isomorphic to $E$, 
such that for any positive and continuous functional $f\colon E\to \RR$ 
there exists a positive and continuous extension $\tl{f}\colon C(X)\ra\RR$. 
Moreover, if $f\in E^{'}_{+}$ 
then $\tl{f}$ can be chosen such that $\tl{f}\in C(X)^{'}_{+}$ and 
$\Vert f \Vert=\Vert \tl{f} \Vert$. \label{isorepn11}

\item There exists a compact Hausdorff space $X$ 
and a subspace $W\subseteq C(X)$, where $W$ is isometrically isomorphic to $E$, 
such that for any functional $f\in E^{'}_{+,1} $ there exists $x_{f}\in K$ such that $\phi(u)(x_{f})=f(u)$ for any $u\in E$. \label{isorepn0}

\item There exists a compact Hausdorff space $X$ 
and a subspace $W\subseteq C(X)$ such that $W$ is 
isometrically isomorphic to $E$. \label{isorepn1}

\item The unit ball 
$E_{1}$ is full; i.e. for $a\leq x \leq b$ with $\Vert a \Vert,\,\Vert b \Vert \leq 1$ we have $\Vert x \Vert\leq 1$. \label{isorepn2}

\item For any $x\in E$ at least one of the following properties holds: $\Vert x+l \Vert \geq \Vert x \Vert$ 
for all $l\in E_+$ and/or $\Vert x-l \Vert \geq \Vert x \Vert$ for all $l\in E_+$, 
i.e. all elements of $E$ are either semi positive and/or semi negative. \label{isorepn3}

\item For all $x,y,z\in E$ with $x\leq y \leq z$ we have $\Vert y \Vert \leq \mathrm{max}\{ \Vert x \Vert, \Vert z \Vert \}$. \label{isorepn9} 

\item For any $u\in E^{'}$, 
there exists $v_1,v_2 \in E^{'}$ such that $u=v_1-v_2$ and $\Vert u \Vert=\Vert v_1 \Vert + \Vert v_2 \Vert$. \label{isorepn4}

\item For any 
$x_0\in E$ there exists a state $f\in E_{+,1}^{'}$ with $\vert f(x_0) \vert = \Vert x_0 \Vert$. \label{isorepn5}

\item For any $x_0\in E$ the equality $\underset{f\in E^{'}_{+,1}}{\mathrm{sup}}\vert f(x_0) \vert = \Vert x_0 \Vert$ holds. \label{isorepn19}

\item There exists a $C^*$-algebra $\cA$ and a real vector subspace $W\subseteq\cA^{h}$ such that $E$ is 
isometrically isometric to $W$. \label{isorepn7}

\item There exists a commutative $C^*$-algebra $\cA$ and a real vector subspace $W\subseteq\cA^{h}$ such that $E$ is 
isometrically isometric to $W$. \label{isorepn18}

\item There exists a Riesz space $(R,\Vert \cdot \Vert)$ with an M-norm $\Vert \cdot \Vert$ 
and a subspace $P\subseteq R$ such that $P$ is isometrically isomorphic to $E$. \label{isorepn10}

\item There exists an almost Archimedean order unit space $(F,F_+,f)$ and a subspace 
$K\subseteq F$ such that $K$ is isometrically isomorphic to $E$. \label{isorepn17}

\item There exists an Archimedean order unit space $(F,F_+,f)$ and a subspace 
$K\subseteq F$ such that $K$ is isometrically isomorphic to $E$. \label{isorepn20}
\end{enumerate}

\end{theorem}

\begin{remark}\label{r:refofisorep}
Most of the items appearing in theorems \ref{t:isorepcx}, \ref{t:isorepsnfl}, \ref{t:isorepstr} and \ref{t:isorepquot} 
are known in the normed case, 
and we state them in Theorem \ref{t:isorep2} as corollaries of the more general calibrated 
ordered space case. However, 
in this article 
the central tools are Schaefer's construction and some consequences of the Bauer-Nachbin-Namioka theorem, 
e.g. see Corollary \ref{c:bau3} in Section \ref{s:phb}  
and generally speaking, the methods we adopt 
in this article are different than the methods of the articles mentioned below. To the best of our knowledge, 
items (\ref{isorep11}) in Theorem \ref{t:isorepcx}, (\ref{isorep51}) and (\ref{isorep19}) in Theorem \ref{t:isorepsnfl} appear 
new characterizations of 
isometrically isomorphic representability of a calibrated ordered space even in the normed case. 
Consequently, items (\ref{isorepn11}), (\ref{isorepn5}) and (\ref{isorepn19}) in Theorem \ref{t:isorep2} 
are new. A similar statement to item $(\ref{isorep60})$ of Theorem \ref{t:isorepstr} was 
given in Theorem 2.13 in \cite{Asadi}; however, as noted in Section \ref{s:nps}, 
the definition of a local order unit space that we give in this article is less restrictive. Consequently, 
the statement obtained in this article is more general. 
 The other 
items are known in the normed case, and in theorems \ref{t:isorepcx}, \ref{t:isorepsnfl}, \ref{t:isorepstr} and \ref{t:isorepquot} 
we state and prove generalizations of them in the locally convex (calibrated) setting. 

We give sources for the known items appearing in Theorem \ref{t:isorep2}. Item (\ref{isorepn0}) is a variation of existence of $C(X)$ with 
some properties. Equivalence of item (\ref{isorepn2}) with isometrically isomorphic representability was
first stated and proved in \cite{Riedl}, see Theorem 4.2. It rests on the Grosberg - Krein result in \cite{GrosbergKrein}, 
which appears as item (\ref{isorepn4}). Item (\ref{isorepn3}) was first stated and its equivalence with isometrically isomorphically 
representability was proved in \cite{Nachbin}. Items (\ref{isorepn7}) and (\ref{isorepn18}) were first observed 
in \cite{Kadison}. For item (\ref{isorepn10}) see \cite{vanGaans}. Items (\ref{isorepn17}) and (\ref{isorepn20}) 
are treated in a number of sources, see e.g. \cite{Kadison} and \cite{PaulsenTomforde}. Finally, 
the definition of a local order unit space was first given in \cite{Dosi} and a corresponding 
representation theorem was obtained in \cite{Asadi}. 
\end{remark} 

We now provide some examples. The following example must be known:

\begin{example}
Any $\ell^p$ space with $1\leq p < \infty$ with the natural cone of positive elements 
$\ell^p_{+}$ of elements with nonnegative entries and the $p$-norm does not 
have a full unit ball, therefore they donot satisfy the equivalent conditions of 
Theorem \ref{t:isorep2} and consequently they donot admit an isometrically 
isomorphic representation on some $C(K)$ with $K$ a compact Hausdorff space. Indeed, for any 
$1 >\epsilon > 0$ consider the elements $a:=(\epsilon^{\frac{1}{p}}, (1-\epsilon)^{\frac{1}{p}}, 0,0,\cdots)$ 
and $b:=(-(1-\epsilon)^{\frac{1}{p}}, -\epsilon^{\frac{1}{p}}, 0,0,\cdots)$. Then both $a$ and $b$ 
are in the unit ball of $\ell^p$. However for $\epsilon$ sufficiently small, the element 
$x:=(-(1/2)^{\frac{1}{p}}, (1/2)^{\frac{1}{p}}, 0,0,\cdots)$ satisfy $a\geq x \geq b$ 
even though $x$ is not in the unit ball of $\ell^p$. On the other hand, these spaces 
have increasing norms, therefore they have normal cones. 

For the case of $\ell^{\infty}$ with the sup norm, it is not difficult to see that it  
has a full unit ball and therefore has all items listed in Theorem \ref{t:isorep2} as
properties. 

An adaptation of these arguments show that, the spaces $L^p(I)$, 
where $I\subseteq \RR$ is an interval, do not have a full unit ball 
for with $1\leq p < \infty$, where the cone consists of almost everywhere nonnegative elements. 
The space $L^{\infty}(I)$, on the other hand, has a full unit ball.
\end{example}

\begin{remark}\label{r:fullext}
For a calibrated ordered space $(E,E_+,\{ p_{\alpha}\}_{\alpha\in A})$ with a closed cone $E_+$ satisfying 
the equivalent conditions of theorems \ref{t:isorepcx}, \ref{t:isorepsnfl}, \ref{t:isorepstr} 
and \ref{t:isorepquot} (for the normed case Theorem \ref{t:isorep2}), in general, a BNN extension of a positive 
and continuous functional does not exist. Namely, given a subspace $V\subseteq E$ and a functional 
$f\colon V\ra \RR$ with $f\in V^{'}_{\alpha, +}$ for some $\alpha$, in general an extension $\tl{f}\colon E\ra \RR$ 
with $\tl{f}\in E^{'}_{\alpha,+}$ does not exist. In \cite{Nachbin} there is an 
example with a space $E$ which does not satisfy the equivalent conditions of Theorem 
\ref{t:isorep2}. Using the same idea, in the following we provide another example satisfying Theorem \ref{t:isorep2}.  
\end{remark}

\begin{example}
Let $E:=\RR^2$, considered as the continuous functions on a two point set with the sup norm and 
a positive cone $E_+:=\{(t_1,t_2) \mid t_1,t_2\geq 0 \}$. Clearly $E_+$ is generating and $E$ satisfies 
the equivalent conditions of Theorem \ref{t:isorep2}. Consider the subspace $X:=\{(t,-2t) \mid t\in \RR  \}$ 
and note that $X_+=E_+\cap X =\{ 0 \}$, so any linear functional on $X$ is positive. Consider 
$f\colon X\to \RR$ by $f(t,-2t):=2t$. It is easy to check that $\Vert f \Vert =1$. In order 
to show that e.g. condition $(3)$ of Theorem \ref{t:baunami} fails on $f/\Vert f \Vert$, 
let $x=(1,-2)\in X$ and $l=(0,1)\in E_+$ and we have $x+l=(1,-1)$, so $x\in X\cap (A-E_+)$ where $A$ is the unit ball. 
But we have $f(1,-2)=2$, so by Theorem \ref{t:baunami} $f$ does not have a BNN extension. 
\end{example}

On the other hand, considering the same algebraic structure equipped with a seminorm, the behaviour 
of BNN extensions may change in the opposite direction. 

\begin{example}
Consider $(E,E_+)$ as in the example above, equipped with the seminorm $p(x,y):=\vert x \vert$. Let 
$X$ and $f$ be as above. Then we have $\vert f(t,-2t) \vert = 2\vert t \vert \leq 2 p(t,-2t)$ and 
it follows that $\vert f \vert_{p}=2$. Let now $(a,t)\in A$, where $A$ is the unit ball, so 
$\vert a \vert\leq 1$ and $t\in \RR$. Let $(t_1,t_2)\in E_+$. 
Then a point $x\in X\cap (A-E_+)$ has the form $(a-t_1, t-t_2)$ with $-2(a-t_1)=t-t_2$. Hence 
$f(a-t_1, t-t_2)=2(a-t_1)$ and $f(x)/\vert f \vert_{p}=a-t_1 \leq 1$. Thus by 
Theorem \ref{t:baunami} $f$ has a BNN extension. Clearly, $f(x,y)=2x$ is such an extension. 
\end{example}

\section{Norm Additivity of Two Positive Functionals}\label{s:naddpos}

In this section we provide an application of the main theorems \ref{t:isorepcx}, \ref{t:isorepsnfl}, \ref{t:isorepstr} and \ref{t:isorepquot}. 
Given a calibrated ordered space 
$(E,E_+,\{ p_{\alpha}\}_{\alpha\in A})$ with a closed cone $E_+$ satisfying 
the equivalent conditions of theorems \ref{t:isorepcx}, \ref{t:isorepsnfl}, \ref{t:isorepstr} and \ref{t:isorepquot}, 
for a fixed $\alpha\in A$, consider two positive functionals 
$f,g\in E^{'}_{\alpha,+}$. Then in general it is not true that 
$\vert f+g \vert_{\alpha}=\vert f \vert_{\alpha}+\vert g \vert_{\alpha}$, 
as can easily be seen by considering the trivial cone $E_+=\{ 0 \}$ and 
nonzero positive functionals $f,-f$. 
In fact, even when $C(X)$ as in $(\ref{isorep1})$ of Theorem \ref{t:isorepcx} 
is considered with an arbitrary generating subcone of the natural cone $C(X)_+$, 
the equality $\vert f+g \vert_{\alpha}=\vert f \vert_{\alpha}+\vert g \vert_{\alpha}$ need not hold, 
as the following example shows. 

\begin{example}
Let $E:=\RR^2$, considered as the continuous functions on a two point set with the sup norm and 
a positive cone $E_+:=\{(t_1,t_2) \mid 4t_2 \leq t_1 \leq 8t_2, \,\,t_1,t_2\geq 0 \}$. It is not difficult to 
check that $E_+$ is generating. Consider the positive functionals 
$g(t_1,t_2):=t_1+t_2$ and $g(t_1,t_2):=t_1-2t_2$. It is not difficult to check that $\Vert f \Vert=2$ 
and $\Vert g \Vert=3$, but $h:=f+g$ is $h(t_1,t_2)=2t_1-t_2$ and consequently $\Vert h \Vert=3$, so 
$\Vert f+g \Vert \neq \Vert f \Vert + \Vert g \Vert$.   
\end{example} 

Below we give a characterization of 
the pairs $f,g\in E^{'}_{\alpha,+}$ 
with the property $\vert f+g \vert_{\alpha}=\vert f \vert_{\alpha}+\vert g \vert_{\alpha}$. 
The space $C(X)$ in the statement and the proof of the theorem refers to the space in 
$(\ref{isorep11})$ of Theorem \ref{t:isorepcx}, note that in this case any positive functional 
$f\in E^{'}_{\alpha,+}$ possess at least one BNN extension $\tl{f}\in C(X)^{'}_{\alpha,+}$, see Section 
\ref{s:phb} for the definition. 

\begin{theorem}\label{t:flsnorm}
Let $(E,E_+,\{ p_{\alpha}\}_{\alpha\in A})$ be a calibrated ordered space with a closed cone $E_+$ satisfying any 
of the equivalent conditions of theorems \ref{t:isorepcx}, \ref{t:isorepsnfl}, \ref{t:isorepstr} and \ref{t:isorepquot}. For a fixed 
$\alpha\in A$, consider two positive functionals 
$f,g\in E^{'}_{\alpha,+}$. Then $\vert f+g \vert_{\alpha}=\vert f \vert_{\alpha}+\vert g \vert_{\alpha}$ 
if and only if the positive functional $f+g$ has a BNN extension $\widetilde{f+g}\in C(X)^{'}_{\alpha,+}$ 
such that $\widetilde{f+g}\geq \tl{f}$ for some BNN extension $\tl{f}\in C(X)^{'}_{\alpha,+}$ of $f$, 
where $C(X)$ is the space in $(\ref{isorep11})$ of Theorem \ref{t:isorepcx}. 
\end{theorem} 

\begin{proof}
``$\Rightarrow$": Assuming $\vert f + g \vert_{\alpha}=\vert f  \vert_{\alpha} + \vert g \vert_{\alpha}$, consider 
the BNN extensions $\tl{f}$ and $\tl{g}$ on $C(X)$, then $\tl{f}, \tl{g} \in C(X)^{'}_{\alpha,+}$. 
It is clear that $\tl{f} + \tl{g}$ is a positive extension of $f+g$, and also clearly $\tl{f} + \tl{g} \geq \tl{f}$. 
We claim that $\tl{f} + \tl{g}$ is a BNN extension of $f+g$. We have 
$\vert \tl{f} \vert_{\alpha} + \vert \tl{g} \vert_{\alpha} \geq \vert \tl{f} + \tl{g} \vert_{\alpha} 
\geq \vert f + g \vert_{\alpha} 
= \vert f  \vert_{\alpha} + \vert g \vert_{\alpha} 
= \vert \tl{f} \vert_{\alpha} + \vert \tl{g} \vert_{\alpha}$ 
where the second inequality follows since $\tl{f} + \tl{g}$ is an extension of $f+g$, and the last equality follows 
since $\tl{f}$ and $\tl{g}$ are BNN extensions of $f$ and $g$, respectively. It follows 
that $\vert \tl{f} + \tl{g} \vert_{\alpha} = \vert f + g \vert_{\alpha}$ and therefore 
$\tl{f} + \tl{g}$ is a BNN extension of $f+g$. 

``$\Leftarrow$": Let us denote $h:=f+g$ and $\tl{h}:=\widetilde{f+g}$, then $\tl{h}\geq \tl{f}$, 
note that $\vert \tl{h} \vert_{\alpha}=\vert h \vert_{\alpha}$ and 
$\vert \tl{f} \vert_{\alpha}=\vert f \vert_{\alpha}$ since 
$\tl{h}$ and $\tl{f}$ are BNN extensions. Consider 
the quotient space $(C(X))_{p_{\alpha}}$ and 
note that by the discussion preceeding Theorem \ref{t:rmlocal}, $(C(X))_{p_{\alpha}}$ isometrically isomorphically identified with 
$C(X_{\alpha})$ for some 
compact Hausdorff space $X_{\alpha}$. Since $E\subseteq C(X)$ isometrically isomorphically by assumption, the quotient 
space $(E_{\alpha}, \ol{(E_+)_{\alpha}}, \Vert \cdot \Vert_{\alpha})$ is 
isometrically isomorphicly identified with a subspace of $C(X_{\alpha})$. Consider the quotient 
functionals $h_{p_{\alpha}}$, $f_{p_{\alpha}}$ and $g_{p_{\alpha}}$ on $(E)_{p_{\alpha}}$ as 
in the discussion preceeding Theorem \ref{t:rmlocal}. Passing to the quotient functionals $\tl{h}_{p_{\alpha}}$ and 
$\tl{f}_{p_{\alpha}}$ on $C(X_{\alpha})$ as well, $\tl{h}_{p_{\alpha}}$ and 
$\tl{f}_{p_{\alpha}}$ are positive functionals extending $h_{p_{\alpha}}$ and $f_{p_{\alpha}}$, 
respectively and by the 
discussion preceeding Theorem \ref{t:rmlocal} we have $\Vert \tl{h}_{p_{\alpha}} \Vert=\vert \tl{h} \vert_{\alpha}$ 
and $\Vert \tl{f}_{p_{\alpha}} \Vert=\vert \tl{f} \vert_{\alpha}$; in particular, 
$\tl{h}_{p_{\alpha}}, \tl{f}_{p_{\alpha}} \in (C(X_{\alpha}))^{'}_{+}$. Also, $\tl{h}_{p_{\alpha}} - \tl{f}_{p_{\alpha}} \geq 0$ 
in $C(X_{\alpha})^{'}_{+}$   
since $\tl{h}_{p_{\alpha}}(x_{p_{\alpha}}) - \tl{f}_{p_{\alpha}}(x_{p_{\alpha}}) = \tl{h}(x) - \tl{f}(x) \geq 0$ for all 
$x\in C(X)_+$ by the assumption. 

Now let   
$\cC(X_{\alpha}):=C(X_{\alpha})+ iC(X_{\alpha})$, the space of all $\CC$-valued continuous functions 
on the compact space $X_{\alpha}$. Then by Section \ref{s:cirpolcs}, $\cC(X_{\alpha})$ is a commutative $C^*$-algebra. 
Let us define a $\CC$-linear functional $h^{\dagger}\colon \cC(X_{\alpha}) \ra \CC$ by $h^{\dagger}:=\tl{h}_{p_{\alpha}} + i \tl{h}_{p_{\alpha}}$, i.e. 
$h^{\dagger}(x_{p_{\alpha}}+ i y_{p_{\alpha}}):= \tl{h}_{p_{\alpha}}(x_{p_{\alpha}}) + i \tl{h}_{p_{\alpha}}(y_{p_{\alpha}})$ 
for any $x,y\in C(X)$.  
The fact that $h^{\dagger}$ is $\CC$-linear is a routine verification and is similar to e.g. the 
proof of Proposition 3.10 in \cite{PaulsenTomforde}. Similarly, 
define $f^{\dagger}\colon \cC(X_{\alpha}) \ra \CC$ by $f^{\dagger}:=\tl{f}_{p_{\alpha}} + i \tl{f}_{p_{\alpha}}$. 
Clearly, $h^{\dagger}$ is an extension of $\tl{h}_{p_{\alpha}}$ and $f^{\dagger}$ is an extension of $\tl{f}_{p_{\alpha}}$. 
We also have $\Vert h^{\dagger} \Vert = \Vert \tl{h}_{p_{\alpha}} \Vert$: Note that, by definition of the operator norm 
we have $\vert \tl{h}_{p_{\alpha}}(x_{p_{\alpha}})\vert \leq \Vert \tl{h}_{p_{\alpha}} \Vert \Vert x_{p_{\alpha}} \Vert$, 
where $\Vert x_{p_{\alpha}} \Vert=\underset{t\in X_{\alpha}}{\mathrm{sup}}\,\vert x_{p_{\alpha}} (t) \vert$ for any 
$x\in C(X)$. It follows that, 
\begin{align*}
\vert h^{\dagger} (x_{p_{\alpha}}+ i y_{p_{\alpha}})\vert &= \vert \tl{h}_{p_{\alpha}}(x_{p_{\alpha}}) + i \tl{h}_{p_{\alpha}}(y_{p_{\alpha}}) \vert \\
&=\sqrt{\vert \tl{h}_{p_{\alpha}}(x_{p_{\alpha}}) \vert^2 + \vert \tl{h}_{p_{\alpha}}(y_{p_{\alpha}}) \vert^2} \\
&\leq \Vert \tl{h}_{p_{\alpha}} \Vert \underset{t\in X_{\alpha}}{\mathrm{sup}}\,\sqrt{\vert x_{p_{\alpha}} (t) \vert^2+ \vert y_{p_{\alpha}} (t) \vert^2} 
\end{align*}
for any $x,y\in C(X)$ and consequently $\Vert h^{\dagger} \Vert \leq \Vert \tl{h}_{p_{\alpha}} \Vert$. Since $h^{\dagger}$ is an extension 
of $\tl{h}_{p_{\alpha}}$, we have $\Vert \tl{h}_{p_{\alpha}} \Vert \leq \Vert h^{\dagger} \Vert $ and thus 
we have $\Vert h^{\dagger} \Vert = \Vert \tl{h}_{p_{\alpha}} \Vert$, 
thus $\Vert h^{\dagger} \Vert = \Vert \tl{h}_{p_{\alpha}} \Vert=\vert \tl{h} \vert_{\alpha}=\vert h \vert_{\alpha}$. A similar argument shows that 
$\Vert f^{\dagger} \Vert = \Vert \tl{f}_{p_{\alpha}} \Vert=\vert \tl{f} \vert_{\alpha}=\vert f \vert_{\alpha}$. 

Since $f^{\dagger}$ and $h^{\dagger} - f^{\dagger}$ are positive functionals on the $C^*$-algebra $\cC(X)$, by 
a well known result from the theory of $C^*$-algebras, we have 
$\Vert h^{\dagger} \Vert = \Vert f^{\dagger} \Vert + \Vert h^{\dagger} - f^{\dagger} \Vert$, see 
Corollary 3.3.5 in \cite{Murphy}. Equivalently we have 
$\vert h \vert_{\alpha}=\vert f \vert_{\alpha} + \Vert h^{\dagger} - f^{\dagger} \Vert$. Since 
$h^{\dagger} - f^{\dagger}$ is an extension of $\tl{h}_{p_{\alpha}} - \tl{f}_{p_{\alpha}}$ 
and $\tl{h}_{p_{\alpha}} - \tl{f}_{p_{\alpha}}$ is an extension of $h_{p_{\alpha}} - f_{p_{\alpha}} = g_{p_{\alpha}}$, 
we have 
\begin{align*}
\vert g \vert_{\alpha} &= \Vert g_{p_{\alpha}} \Vert = \Vert h_{p_{\alpha}} - f_{p_{\alpha}} \Vert \\
& \leq \Vert \tl{h}_{p_{\alpha}} - \tl{f}_{p_{\alpha}} \Vert \\
& \leq \Vert h^{\dagger} - f^{\dagger} \Vert
\end{align*}
from which we obtain $\vert g \vert_{\alpha} + \vert f \vert_{\alpha} \leq \vert h \vert_{\alpha}$. 
Since the reverse inequality is clear, we obtain $\vert g \vert_{\alpha} + \vert f \vert_{\alpha} = \vert h \vert_{\alpha}$, 
as required.  
\end{proof}

An inspection of the proof of Theorem \ref{t:flsnorm} shows that, the if part 
of the theorem holds with BNN extensions considered not necessarily in the space 
as in $(\ref{isorep11})$ of Theorem \ref{t:isorepcx}, but in an arbitrary space 
$C(X)$ with $X$ a locally compact Hausdorff space. Consequently, 
we obtain the following corollary:

\begin{corollary}\label{c:flsnorm}
Let $E$ be an ordered subspace, with a closed cone, of the calibrated ordered space $(C(X),C(X)_+,\{ K_{\alpha}\}_{\alpha\in A})$ where 
$X$ a locally compact Hausdorff space. For a fixed 
$\alpha\in A$, consider two positive functionals 
$f,g\in E^{'}_{\alpha,+}$ and assume that a BNN extension $\tl{f}\colon C(X)\to \RR$ of 
$f$ exists. Then $\vert f+g \vert_{\alpha}=\vert f \vert_{\alpha}+\vert g \vert_{\alpha}$ 
if the positive functional $f+g$ has a BNN extension $\widetilde{f+g}\in C(X)^{'}_{\alpha,+}$ 
such that $\widetilde{f+g}\geq \tl{f}$.
\end{corollary}

\begin{remark}
In the case when a calibrated ordered space $E$ is isometrically isomorphic with 
$C(X)$, where $X$ is locally compact and Hausdorff, see Remark \ref{r:onto}, the condition 
of Corollary \ref{c:flsnorm} is trivially satisfied, therefore for such spaces 
seminorm additivity on positive functionals holds automatically. In particular, it is 
well known that AM (Abstract M) spaces with unit are such spaces, see \cite{Kakutani}, 
as well as \cite{Schaefer1} and \cite{Riedl}. 
\end{remark}


\end{document}